\documentclass[11pt, leqno]{article}
\usepackage{amsmath,amsthm}
\usepackage{amsfonts}
\usepackage{mathrsfs}
\usepackage{amssymb}
\usepackage{enumitem}
\usepackage{epsfig}
\usepackage{esint}    
\usepackage{pst-grad} 
\usepackage{pst-plot} 
\usepackage{pinlabel}
\usepackage{graphicx}
\usepackage{color}

\newtheorem{theor}{Theorem}[section]
\newtheorem{cor}[theor]{Corollary}
\newtheorem{lemma}[theor]{Lemma}
\newtheorem{prop}[theor]{Proposition}

\newtheorem*{question*}{Question}

\theoremstyle{definition}

\theoremstyle{remark}

\numberwithin{equation}{section}
\numberwithin{defn}{section}

\newcommand{\R}{\mathbb{R}}        
\newcommand{\re}{\mathbb{R}}        
\newcommand{\RN}{\mathbb{R}^N}
\newcommand{\ren}{\mathbb{R}^N}

\newcommand{\eps}{\varepsilon}

\newcommand{\ve}{\varepsilon}

\DeclareMathOperator{\Lap}{(-\Delta)} 

\newcommand{\Ds}{\Lap^{s}}

\newcommand{\Doh}{\Lap^{1/2}}
\renewcommand{\d}{\mathrm{d}}

\definecolor{darkblue}{rgb}{0.05, .05, .65}
\definecolor{darkgreen}{rgb}{0, .6, .2}
\definecolor{darkred}{rgb}{0.8,0,0}

\begin{document}
\pagenumbering{roman}
\title{ \bf Existence of maximal solutions \\for some very singular nonlinear  \\ fractional diffusion equations in 1D}
\author{ \ Juan Luis V\'azquez \\
Dpto. de Matem\'aticas \\ Univ. Aut\'onoma de Madrid}

\date{}

\maketitle

\begin{abstract}
We consider nonlinear parabolic equations involving fractional diffusion of the form $\partial_t u + (-\Delta)^s \Phi(u)= 0,$ with $0<s<1$,  and solve an open problem concerning the existence of solutions for very singular nonlinearities $\Phi$ in power form, precisely $\Phi'(u)=c\,u^{-(n+1)}$ for some $0< n<1$. We also include the logarithmic diffusion equation $\partial_t u + (-\Delta)^s \log(u)= 0$, which appears as the case $n=0$. We consider the Cauchy problem with nonnegative and integrable data $u_0(x)$ in one space dimension, since the same problem in higher dimensions admits no nontrivial solutions according to recent results of the author and collaborators. The {\sl limit solutions} we construct are unique, conserve mass, and are in fact maximal solutions of the problem. We also construct self-similar solutions of Barenblatt type, that are used as a cornerstone in the existence theory, and we  prove that they are asymptotic attractors (as $t\to\infty$) of the solutions with general integrable data. A new comparison principle is introduced.

\end{abstract}

\


\medskip


\vskip.5cm

\small
\setcounter{tocdepth}{1}
\tableofcontents
\normalsize

\pagenumbering{arabic}
\setcounter{page}{1}
\section{Introduction}

In this paper we consider a class of nonlinear parabolic equations involving fractional diffusion  of the form
\begin{equation}\label{geq.Phi}
\partial_t u + (-\Delta)^s \Phi(u)= 0\,.
\end{equation}
The symbol  $\Ds$ denotes  the fractional Laplacian operator with $0<s<1$, i.\,e.,  the nonlocal operator defined  by
\begin{equation}
\label{eq:fract_lap}
(-\Delta)^s v (x) = c(N,s)\, \hbox{p.v.}\int_{\RN}\frac{v(x) - v(y)}{\vert x-y\vert^{N+2s}}\d y, \,\,\,\forall x\in \RN\,,
\end{equation}
acting on the whole Euclidean space at least for functions in the Schwartz class $\mathcal{S}$. The formula is valid for all dimensions $N\ge 1$. The constant $c(N,s)$ is given in the literature but it is not needed in what follows and p.v. means principal value of the integral.

The existence and properties of solutions for this type of equations with fractional diffusion has  been studied by the author and collaborators for nonlinearities $\Phi$ that are positive and increasing for $u>0$, in particular when $\Phi(u)=u^m$ with $m>0$, cf. \cite{pqrv1, pqrv2, pqrv3,  BV2012, VazVol1, VPQR13}. This includes singular cases for $0<m<1$ since then $\Phi'(u)=mu^{m-1}\to\infty$ as $u\to 0$.

Here, we are interested in {\sl very singular nonlinearities}, more precisely,  when $\Phi:\R_+\to \R$ is a monotone increasing function of $u$ with a singularity in $u=0$ such that $\Phi(0+)= -\infty$. Consequently, nonnegative data and solutions are considered. The standard cases we have in mind are \ $\Phi_{n}(u)=-1/u^n$ for some $n>0$, or $\Phi_0(u)=\log(u)$. They correspond to $\Phi'(u)\sim u^{-(n+1)}$ with $n+1\ge 1$, thus the denomination {\sl very singular} introduced in the literature for this type of equations with standard Laplacian, cf. \cite{JLVSmoothing}. We will keep this tradition for equations with a fractional Laplacian. These very singular diffusion equations are also described in the literature as {\sl very fast diffusion},  {\sl superfast diffusion}, or {\sl ultra-fast diffusion}, cf. e.\,g. \cite{Igb09, Rosen95, JLVSmoothing}.

For such equations the existence of solutions is not at all obvious. Thus, we have proved  in a recent paper with Bonforte and Segatti \cite{BSV15} that when the space dimension is $N\ge2$ and we try to solve the Cauchy problem in the whole space $\ren$ with integrable initial data, then there exist no nontrivial solutions, even if we accept local-in-time solutions defined for a short time interval, $0<t<T$.\footnote{ By trivial solution we mean $u(x,t)\equiv 0$ in the whole domain of definition.} The same happens for the problem posed in a bounded domain with zero Dirichlet data.

The purpose of this paper is to prove that there is a range of existence of solutions for very singular equations of the form \eqref{geq.Phi} with $\Phi=\Phi_n$ if the space dimension is 1.  We will also prove that the solutions have the good properties of the non-singular range of parameters $\Phi(u)=u^m$ with $m>0$. A very crude explanation of the existence result is as follows: by becoming strictly positive for $t>0$, the solutions avoid the singularity in a way that suffices to grant first nontrivial existence, and then the rest of the properties.

Before stating the results, let us point out that the standard notation for the nonsingular equation is $\Phi(u)=c\,u^m$ with $c, m>0$. In this paper the exponent of the nonlinearity $\Phi_n$ is written in terms of $n=-m$ with $n\ge 0$, and  the detailed calculations are done  for $n>0$. The reason for this notation is to avoid the use of negative exponents that might confuse the reader in interpreting the results\footnote{The reader may also wonder, why the minus sign in the coefficient of $\Phi_n$? It is needed to make $\Phi_n$ an increasing function, so that the equation will be parabolic in some sense.}.  The space dimension is $N$, mostly $N=1$ here.

\begin{theor} \label{thm.1} Let $N=1$ and let $\Phi(u)=\Phi_{n}(u)$ with $n> 0$ or $\Phi(u)=\log(u)$ (case $n=0$). Equation  \eqref{geq.Phi} posed in $Q=\re\times (0,\infty)$ with initial data
\begin{equation}\label{id}
u(x,0)=u_0(x)\in L^1(\re), \   \quad u_0\ge 0\,,
\end{equation}
admits a positive very weak solution if \ $ s>1/2$ and $0\le n < 2s-1$.
There could be non-uniqueness of the solutions, but we construct a unique limit solution  for every initial data, and we prove that it is maximal among all solutions. This solution exists globally in time, $u\in C([0,\infty):L^1(\re))$, and is positive everywhere.
\end{theor}

The range of exponents $ 1/2<s<1$ and $0\le n < 2s-1$ is almost sharp. There is indeed another isolated case of existence of integrable solutions in 1D, namely $s=1/2$ and $n=0$ (logarithmic diffusion). The very peculiar properties of this case deserve a separate study, but we give a preliminary idea in Section  \ref{sec.logdiff} that supports the assertion of existence of solutions at least for short times. For the other exponents $s\in (0,1)$, $n\ge 0$, in dimension one, nontrivial solutions do not exist by the mentioned results of \cite{BSV15}.

The construction of solutions proceeds by approximation, taking approximate initial data that are uniformly positive, so that the problem is no more singular. Passing then to the limit in the approximations we obtain a solution that is shown to be non-trivial after some effort. It is called the {\sl limit solution} ({\sl upper limit solution}, to be precise). It is subsequently proved to be a very weak solution. Here we define very weak solution of equation \eqref{geq.Phi} as a nonnegative function $u\in C((0,\infty):L^1(\ren))$ such that
\begin{equation}\label{weak-nonlocal}
\displaystyle \int_0^\infty\int_{\re}u\dfrac{\partial \zeta}{\partial t}\,dxdt
=\int_0^\infty\int_{\re}  \Phi(u)\,(-\Delta)^{s}\zeta\,dxdt\,,
\end{equation}
and the last integral is absolutely convergent for all $\zeta$ smooth and compactly supported.
In the theorem $\Phi(u)=\Phi_n(u)$.

In the course of the present paper we will also establish the main properties of the constructed solutions. We select here the main results for easy reference.

\begin{theor} \label{thm.2} The limit solution  preserves mass, $\int u(x,t)\,dx=\int u_0(x)\,dx$. Moreover, it is  a weak solution for $t>0$, and satisfies the bounds
 \begin{equation}\label{bounds}
 C_1(t)(1+|x|^2)^{s/(1+n)}\le u(x,t)\le C_2\|u_0\|_1^\delta t^{-\alpha}
 \end{equation}
with $\alpha=1/(2s-n-1)$, $\delta=2s\alpha$, a continuous function $C_1(t)> 0$ that may depend on the solution, and a constant $C_2(n,s)>0.$ The collection of limit solutions generates an ordered, $L^1$-contraction semigroup in $L^1_+(\re)$
\end{theor}

 A specially important feature of the paper is the construction and properties of the fundamental solutions.

\begin{theor} \label{thm.3} There exists a special function of the form
\begin{equation}
U(x,t)=t^{-\alpha} F(x\,t^{-\alpha})
\end{equation}
that  is a very weak positive solution of the problem for $t\ge \tau>0$ and takes on a Dirac mass as initial data,  $u(x,t)\to \delta(x)$ as $x\to 0$ in the sense of positive Radon measures.
The profile $F$ is positive everywhere, integrable, symmetric, $F(x)=F(-x)$, and $F$ monotone decreasing for $x>0$.  Moreover,
\begin{equation}
\lim_{|x|\to \infty } |x|^{2s/(1+n)} F(|x|)= C(s,n)>0.
\end{equation}
 The constant $C(s,n)$ can be calculated as the constant appearing in the Very Singular Solution, a special solution with formula $\widetilde U(x,t)=C(s,n)\,t^{1/(1+n)}|x|^{-2s/(1+n)}$, that has a non-integrable singularity at $x=0$.
 Finally, the solution with initial data $M\,\delta (x)$, $M>0$, is just
 \begin{equation}
  U_M(x,t)=M\,U(x,M^{-(1+n)}t)\,, \quad \mbox{so that } \ F_M(\xi)=M^{2s\alpha}\,F(M^{(1+n)\alpha}\xi).
 \end{equation}
 \end{theor}

Using the terminology of \cite{VazBar} we call a {\sl fundamental} any solution with $u(x,0)$ equal to a Dirac delta, and Barenblatt solution a fundamental solution that is also self-similar. On the one hand, the Barenblatt solutions of the theorem play an important role in completing the existence theory described in Theorem \ref{thm.2}. In the theory different comparison theorems are also used, in particular a new Shifting Comparison result, that we prove as Theorem \ref{thm.shift}. On the other hand, the Barenblatt solutions explain the asymptotic behaviour of general solutions, according to the following general theorem.

\begin{theor}\label{thm.4} Let $u_0\in L^1(\re)$, let $M=\int u_0(x)\,dx$, and let $U_M$ be the self-similar Barenblatt solution with mass $M$. Then  as $t\to\infty$ the solutions $u(x,t)$ and $U_M(x,t)$ are increasingly close and we have
\begin{equation}\label{conver.express.1}
\lim_{t\to\infty} \|u(\cdot,t)-U_M(\cdot,t)\|_1=0\,,
\end{equation}
Indeed, convergence happens in all $L^p$ norms, $1\le p<\infty$, in the form
\begin{equation}
\lim_{t\to\infty} t^{\alpha_p}\|u(\cdot,t)-U_M(\cdot,t)\|_{L^p(\re)}=0\,, \qquad
\alpha_p=   \frac{p-1}{p(2s-1+n)}.
\end{equation}
\end{theor}

There are some other results worth recalling. Thus, as a side result of our analysis, we construct in Section \ref{sec.vss} the Very Singular Solution (VSS), that is explicit (see Theorem \ref{thm.vss2}) and has very special properties. Very Singular Solutions have played a special role in the theory of fast diffusion equations, as attested e.\,g. in \cite{ChaVaz2002}. Our VSS will give us a first clue to the lower bound $O(|x|^{-2s/(1+n)})$ for the spatial decay of all positive solutions, that we have stated in \eqref{bounds} and plays a role in different passages of the existence theory that we will develop below.

After all these theorems are proved and shifting comparison is established, we devote a short section to a preliminary presentation of the special case $s=1/2$, $n=0$. The paper concludes with a section on comments, extensions, and open problems.

\medskip

\noindent {\bf Precedents and commentary.} (1) Many results are known about Problem \eqref{geq.Phi}-\eqref{id}, mainly for standard diffusion $s=1$, where  the Laplacian is used instead of the fractional Laplacian. The nonsigular case $\Phi(u)=u^m$ with $m>0$, is known as the Porous Medium Equation when $m>1$, the Heat Equation for $m=1$ and the Fast Diffusion Equation, $0<m<1$; their theory has been studied in great detail and is described in monographs like \cite{ArBk86, DaskKenig, Vapme, JLVSmoothing}. As a basic existence result, each of these equations generate a mild solution for every initial data $u_0\in L^1(\ren)$ and the collection of such solutions forms an ordered $L^1$ contraction semigroup for every fixed $m$.

\noindent (2) For equations with fractional Laplacians of the form \eqref{geq.Phi}-\eqref{eq:fract_lap} with $0<s<1$ and the same of power-like nonlinearity $\Phi(u)=u^m$, $m>0$, the study of the Cauchy problem with nonnegative data in $L^1(\ren)$, $N\ge 1$, has been done in the papers \cite{pqrv1, pqrv2}, and most of the basic results are still true though the techniques may be quite different. More precisely, existence and uniqueness of solutions in the class of very weak or strong solutions have been proved, and the main qualitative and quantitative properties are established. Thus,  when $N(m-1) + 2s>0$  the solutions are positive everywhere in $Q=\ren\times (0,\infty)$; the so-called smoothing effect asserts that $L^1$ initial data produce bounded solutions for positive times and indeed
$$
u(x,t)\le C(N,m,s)\|u_0\|_1^\gamma\,t^{-\alpha }
$$
where $\alpha=N/(N(m-1)+2s)$ and $\gamma= 2s/(N(m-1)+2s)$, both positive in this range. A main feature of the theory is the existence of fundamental solutions and their use in establishing the asymptotic behaviour of general solutions. This was proved in \cite{VazBar}, also under the necessary restriction $N(m-1) + 2s>0$.

\noindent (3) Since this condition on the exponents to obtain good behaviour becomes  formally $m>1- 2s$ in 1D,  a similar theory could be expected to hold for very singular exponents $m<0$, if $2s-1>n=-m>0$ when $N=1$. However, the difficulties of dealing the singular nonlinearities prevented the inclusion of this extension in the works \cite{pqrv1, pqrv2} and \cite{VazBar}. We supply in this paper the approach and tools to fill such a gap in the range of exponents of Theorem \ref{thm.1}. In particular, we construct the fundamental solutions and show that they are responsible for the asymptotic behavior of general solutions. The text below shows that such extended theory is far from immediate and needs some involved tools.

\noindent (4)  The case of singular powers $\Phi(u)=-u^{-n}$ with $n>0$, or $\Phi_0(u)=\log(u)$, was considered by the author many years ago in \cite{Vaz92} for the standard Laplacian, $s=1$.  A remarkable result of non-existence of integrable solutions was proved for all $n \ge 0$ if $N\ge 3$, for $n>0$ if $N=2$, and for $n\ge 1$ if $N=1$. That paper is a remote precedent for the present work (previously, non-existence for the particular limit case $n=1$ in $N=1$ had been proved in \cite{Herrero89} using a special transformation). More precisely, when we perform the natural approximation by regular problems, as explained in the next section, the sequence of approximations collapses to zero for all $x\in\re$ and all $t>0$ for all initial data in the integrable class. This radical phenomenon is  called instantaneous extinction. On the other hand, existence occurs in the remaining cases, $N=1$, $0\le n<1$, and $N=2$, $n=0$. This agrees with the results we are going to prove in the fractional case. Interesting properties arise in the existence cases, see a detailed account in \cite[Section 9]{JLVSmoothing}.  The non-existence results of \cite{Vaz92} were extended to optimal classes of (non-integrable) initial data in  \cite{DDP1994, DDP1995, DDP1997}.

\noindent (5) In the case of singular nonlinearities and fractional Laplacians, which is our framework here, the non-existence of solutions has been recently established  in collaboration with Bonforte and Segatti,  \cite{BSV15},  in the range that perfectly complements the positive result of Theorem \ref{thm.1} plus the announced existence result for $s=1/2$, $n=0$ in 1D. Non-existence happens for all singular power cases ($n\ge 0$) in dimension $N\ge 2$. In all those non-existence cases we got {\sl instantaneous extinction} for the initial value problem with any data  $u_0\in L^1_+(\re)$. Putting these results together, we obtain a complete picture of the solubility problem with integrable data for all singular parameters.

\medskip

\noindent {\sc More on notations.} We use the sign $f \sim g$ to denote that both functions are proportional in a certain limit or range of values $I$ (which may be explicit or understood from the context). If the proportionality ratio goes to 1 in some limit then we write $f\approx g$.
In the proofs we will often use {\sl rearranged} functions defined on the line. This means that they are nonnegative, symmetric and monotone nonincreasing for $x>0$. Other notations will be explained as they appear.

\section{Problem, approximation, and limit solutions}\label{sec.limit}

Let us discuss the way to prove existence for the  Cauchy Problem \eqref{geq.Phi}--\eqref{id}, i.\,e., to find solutions of the equation posed in $Q=\re\times(0,\infty)$ with $T>0$,  taking on initial data $u_0(x)$, assumed to be nonnegative and integrable. In this section we also assume that $u_0$ is bounded, a restriction that is made for convenience and will be removed later on. The nonlinear function $\Phi$ is defined, increasing and smooth for $u>0$, with $\Phi(s)\to -\infty$ as $s\to0$. More precisely, we will construct a compete theory for the case where $\Phi$  is chosen from the list $\Phi_{n}$, $n\ge 0$, mentioned in the Introduction. From this moment on we assume $n$ to be fixed.
Following \cite{Vaz92}, a strategy of proof of existence or non-existence of solutions is based on approximating problem \ \eqref{geq.Phi}--\eqref{id}  by the family problems
\begin{equation}
\label{eq:FrUF_approx}
\begin{cases}
\partial_t u_\ve  + (-\Delta)^s ({\Phi(u_\ve)}) = 0, \qquad n>0,\\
u_\ve(x, 0) = u_0(x) + \ve\,\,\, \hbox{ for } x\in \RN,
\end{cases}
\end{equation}
for any $\ve>0$, so that we avoid data with values on the singular level $u=0$. The standard theory  applies to these problems and a classical solution $u_\ve(t,x)$  exists for all $\ve>0$, and it is strictly positive: $u_\ve\ge \ve$ (see details in the next subsection). Moreover, the maximum principle holds for these classical solutions and we have $u_\ve\ge u_{\ve'} $ for $\ve\ge \ve'>0$. Therefore, we can take the monotone limit
\begin{equation}\label{lim.sol}
\bar u(x,t)=\lim_{\ve\to 0} u_\ve(x, t)\,.
\end{equation}
This function is a kind of generalized solution of the problem, that belongs to the class of {\sl limit solutions.} It is now an important step of the theory to decide in which  sense this limit solution is a solution of the equation in a more traditional functional sense (like very weak, weak, strong or viscosity solution), and also in which sense it takes the initial data. In cases of non-uniqueness of such solutions, the unique limit obtained by the above method has been called by various names: maximal solution, SOLA, proper solution,... However, this kind of considerations are not the main issue of this paper which is concerned with describing the existence and behaviour of the class of limit solutions.
In order to recall the way the limit is taken  (via approximations from above) and to avoid possible confusions, we propose the more precise term {\sl  upper limit solutions} for the limits \eqref{lim.sol}, but we will allow the simpler name limit solutions when there is no fear of confusion.

This approximation method has been used in  \cite{Vaz92} to prove non-existence in the case of standard Laplacian as mentioned above, and in \cite{BSV15} to prove the non-existence results for fractional diffusion and those singular $\Phi$ that do not fall into the cases treated in this paper.


\subsection{Existence and properties of the approximate solution}

It is convenient to write $u_\ve(x,t)=v_\ve(x,t) +\ve$ and then try to solve
 the Cauchy problem
\begin{equation}
\label{eq:FrUF_approx1}
\begin{cases}
\partial_t v + (-\Delta)^s (\Phi_\ve(v)) = 0\,\,\,\,\,\hbox{ with }\Phi_\ve(v):= \Phi(v+\ve) -\Phi(\ve)
\\
v(0) = u_0 \,\,\,\,\hbox{ for } x\in \R.
\end{cases}
\end{equation}
for all $\ve>0$. This is a modified problem prepared to avoid the singular level $u=0$ of the equation by displacement of the axes. Note that for $\ve>0$  and for nonnegative arguments $\Phi_\ve$ is a smooth, positive, monotone increasing function with $\Phi_\ve(0) = 0$ and $\Phi'_\ve(v)$ positive, bounded and decreasing for all $v\ge 0$; $\Phi'_\ve(v)$ is uniformly positive if $v$ is bounded. The theory of existence and uniqueness of weak solutions  $v_\ve$ to the nonsingular Problem \eqref{eq:FrUF_approx1} is given in \cite[Theorem 8.2]{VPQR13}.  After obtaining these solutions  we restore  for any $\ve>0$ the original $u$-level by defining
$u_\ve:=v_\ve +\ve$, as stated at the beginning.
Clearly, we have that $u_\ve\ge \ve$ (actually, $u_\ve> \ve$) and
hence $v_\ve\ge 0$ in $Q$.

Let us list some further properties of $v_\ve$ and $u_\ve$. For the proof we may again refer to  \cite{pqrv2,VPQR13} and \cite{BSV15}.

\noindent $\bullet$ {\sl Boundedness and regularity.} These  solutions are shown to be bounded for strictly positive times. More precisely, for every $t>0$ and every $p\in [1,\infty]$ there holds
\begin{equation} \label{eq:Lp}
\|v_\ve(\cdot,t)\|_p\le \|u_0\|_p
\end{equation}
As a consequence,  if $u_0$ belongs to $L^\infty$, then $v_\ve$ is regular enough to satisfy the equation in the classical sense at least when $t>0$\,, by the results of \cite{VPQR13}. Therefore, $u_\ve=v_\varepsilon+\varepsilon$ is smooth and satisfies the original equation in the classical sense in $ Q$. Under these circumstances, the initial data are also taken, at least in the sense of convergence in $L^1(\RN)$. For unbounded data this result about initial data follows from density and contraction in $L^1$, see next paragraph.

\noindent $\bullet$ {\sl $L^1$-contraction and comparison. } The evolution \eqref{eq:FrUF_approx1} is an $L^1$ contraction, namely for two solutions $u_{1,\ve}, u_{2,\ve}$ we have
we have
\begin{equation}\label{eq:contraction}
\int_{\RN}(u_{1\ve}(x,t) -u_{2\ve}(x,t))_{+}\d x \le\int_{\RN}(u_{01} -u_{02})_{+}\d x \,\,\,\hbox{ for }t>0,
\end{equation}
Here,  $(\cdot)_+$ denotes the positive part function. In particular, standard comparison follows: if $u_{01}\le u_{02}$ a.e., then for every $t>0$ we get $u_{1\ve}(\cdot, t)\le u_{2\ve}(  \cdot, t)$ a.e.

\noindent $\bullet$ {\sl Mass conservation. } Nonnegative solutions to the evolution equation \eqref{eq:FrUF_approx1} conserve  the mass, cf. \cite{pqrv2,VPQR13}. More precisely, we have for all $t\ge 0$
\begin{equation}
\label{eq:mass_cons}
\int_{\RN} v_\ve(x,t)\,\d x = \int_{\RN} u_0(x)\,\d x\quad \mbox{i.\,e.,}\quad
\int_{\RN} (u_\ve(x,t)- \ve)\,\d x = \int_{\RN} u_0(x)\,\d x\,.
\end{equation}

\noindent $\bullet$ {\sl Monotonicity with respect to \ $\ve$. } An easy version of the above comparison argument shows also that for $0<\ve<\ve'$ we have $0<\ve\le u_\ve\le  u_{\ve'}$.

\noindent $\bullet$ {\sl Time monotonicity.} There is an important monotonicity
property valid of all nonnegative solutions, known as the B\'enilan-Crandall inequality
\begin{equation}\label{eq:AB}
\partial_t u_\ve \le \frac{u_\ve}{(n+1) t}\,\,\,\,\,\forall (x,t)\in Q.
\end{equation}
The argument only uses the scaling invariance of the equation and the maximum principle, so \cite{BenCr}'s argument applies.

\noindent $\bullet$  {\sl The smoothing effect.} It says that all solutions with integrable data are in fact bounded for positive times. This follows from Theorem 8.2 of \cite{VPQR13}, that we adapt to our dimension and notations as follows.

\begin{prop} \label{th:smoothing2}
Let $\Phi\in C^1(\mathbb{R})$ be such that
 $\Phi'(u)\ge C |u|^{-n-1}$ for
some  $n\in\mathbb{R}$  and $|u|\ge C$. If $u_0\in L^1({\mathbb{R}^N})\cap
L^p({\mathbb{R}^N})$, where $p\ge1$ satisfies $2sp>1+n$,  then a weak $L^1$-energy solution to the Cauchy problem~\eqref{geq.Phi}-\eqref{id}  is bounded in $\mathbb{R}\times(\tau,\infty)$ for all  $\tau>0$. More precisely, it satisfies
\begin{equation}
\sup_{x\in{\mathbb{R}^N}}|u(x,t)|\le
\max\{C, \, C_1\,t^{-\alpha_p }\|u_0\|_{p}^{\delta_p}\}
\label{eq:L-inf-p2}\end{equation}
with
$\alpha_p=1/(2sp-n-1)$ \ and \ $\delta_p=2sp\gamma_p$,
the constant $C_1$ depending on $n,\,p,\,\sigma, C$.
\end{prop}

See also  \cite{pqrv2} in that respect. The statement of \cite{VPQR13} does not assume that $m=-n>0$, only that $2sp >1-m$. The assumptions on $\Phi$ are satisfied by the nonlinearities $\Phi_\ve$ of the approximate problems \eqref{eq:FrUF_approx1}, hence the result applies to the approximations $v_\ve$. Recall that we write $m=-n$ and note that the constants in the formula may also depend on $\ve$. We will solve the latter difficulty later on.

\subsection{Passing to the limit} We may now pass to the limit $\ve\to0$ using the monotonicity in $\eps$ of the family $u_\ve$. By the monotone convergence theorem the limit $\bar u$ is taken in the local $L^1$ sense, i.e., in  $L^1(B)$ for every compact subset $B$ of $\RN\times[0,\infty]$. We have

\begin{prop}\label{prop.limit} If $u_0$ is a nonnegative function in $L^1(\RN)$ there exists the monotone limit
$\bar u=\lim_{\ve\to 0} u_\ve$ with local convergence in $L^1(Q)$.
\end{prop}


\section{Existence of nontrivial limit solutions}\label{sec.bddsub}

The main problem with this  procedure concerns the possibility that the limit may become identically zero in $Q_T$ for some relevant class of initial data. This is what happens for $N\ge 2$, as established in \cite{BSV15}. In the cases we study here this failure of existence will not happen as we will show below. In fact, the limit will be nontrivial for all nontrivial initial data. The proof of this general result is long and proceeds in steps. The first of such steps consists in exhibiting at least one nontrivial solution.

A reminder: in the next sections we concentrate on the case of  exponents $s>1/2$ and $0\le n<2s-1$. Note that even if the diffusion is singular, the range is formally the same as the good fast diffusion range $1>m> (N-2s)/N$, considered in the general theory of \cite{pqrv2}, hence it is supercritical in the notation of that paper. But there only exponents $m>0$ were considered.

\medskip

\subsection{The very singular solution, I}

In paper \cite{VazBar} a formal solution of equation $u_t+(-\Delta)^s u^m=0$ is constructed with the form
\begin{equation}
U(x,t)=H(t)F(x)=C(N,s,m)\,t^{1/(1-m)}|x|^{-2s/(1-m)}
\end{equation}
 when $m>(N-2s)/N$, so that the spatial profile has a non-integrable singularity at $x=0$. This type is called {\sl very singular solution} in the literature (VSS for short). It is also proved that such formal solution is a limit of a monotone increasing sequence of standard solutions.

Following that paper we try the same formula here for $m=-n\le 0$ with factors
 $$
 H_1(t)=C\, t^{1/(1+n)}, \quad F_1(x)=|x|^{-2s/(1+n)}.
 $$
\noindent $\bullet$ Let us check that it works in the  parameter range, $s>1/2$,  $n+1<2s$ with $n>0$. We have $\phi(F(x))=-|x|^{2sn/(1+n)}$, so that a simple calculation gives
\begin{equation}
{\cal L}^s\phi(F(x))=-K(s,n)\,|x|^{2sn/(1+n)-2s}= -K(s,n)F(x)\,.
\end{equation}
We prove in the note below that $K(s,n)$ is positive in this range of values of $s, n$. It is now easy to see that if we put
\begin{equation}
U_1(x,t)=H_1(t)F_1(x)\qquad \mbox{ with } \ C^{1+n}=K(n+1),
\end{equation}
then $U_1$ is a solution of the 1D equation in the range $0<n<2s-1$, unless at the singular point $x=0$.

The singularity prevents the VSS from being acceptable as an example of nontrivial solution in the theory we are considering for Problem \eqref{geq.Phi}-\eqref{id}, at this stage. But we will return to this topic in style in Section \ref{sec.vss}.

\smallskip

\noindent {\bf Note on fractional Laplacians of power functions.} The  $s$-Laplacian of a  power, $(-\Delta)^s |x|^{\alpha}$, $\alpha>0$, $\alpha\ne 2s$, is the power $k(\alpha,s)\,|x|^{\alpha-2s},$ with a constant factor that for  $N=1$ equals
\begin{equation}\label{kalpha}
k(\alpha,s)=2^{2s}\frac{\Gamma((1+\alpha)/2)\,\Gamma((-\alpha+2s)/2)}{\Gamma((1+\alpha-2s)/2)\,
\Gamma(-\alpha/2)}
\end{equation}
In our case we take $\alpha=2sn/(1+n)>0$. Hence, $\Gamma((1+\alpha)/2)>0$; besides, $(-\alpha+2s)/2=s/(1+n)>0$ so that $\Gamma((-\alpha+2s)/2) >0$; moreover, $\alpha/2=sn/(1+n)<1$ (since $s<2$ and $n/(1+n)<1/2$), so that $\Gamma(-\alpha/2)<0$. Finally, for our choice of $\alpha$
$$
1+\alpha-2s=\frac{2ns}{1+n}+1-2s=\frac{2sn-(1+n)(2s-1)}{1+n}=\frac{1+n-2s}{1+n}\in (0,1)\,.
$$
It means that $\Gamma((-\alpha+2s)/2)<0$. Therefore, $k(\alpha,s)>0$ for theses particular values of $\alpha$ and $s$. This is what we have called $K(s,n)$ some lines above.

\medskip

\subsection{A bounded subsolution} Though the singularity prevents the VSS from being directly  useful as an example of nontrivial solution,  the idea is not completely lost. We will start from it to construct a useful smooth variant, but it will be only a subsolution.

\begin{lemma}\label{lemma.subsol}  There exists a bounded and smooth  function in separate-variables form
\begin{equation}
\widetilde U(x,t)=H(t)F_2(x)
\end{equation}
which is a positive subsolution of the equation for some $0<t<T$ and all $x\in \R$.
Moreover, $F_2$ is symmetric, radially decreasing and $F_2(x)\approx c |x|^{2s/(1+n)}$ as $x\to \infty$.
Changing the form of $H$ we can get a supersolution. In both cases $H$ is continuous and the initial value $H(0)>0$ can be chosen.\end{lemma}

\noindent

\noindent {\sl Proof.} (i) In order to avoid the problem with the singularity of the VSS, we round  the function $F_1$ in a small ball near $x=0$ to get a smooth positive $F_2$, so that instead of the exact formula $(-\Delta)^s\phi(F_1)=-K(s)F_1(x)$ for $x\ne 0$ we get an approximate equation for all $x$ that we can still use. Indeed, in view of the perturbation   we get that
$(-\Delta)^s\phi(F_2)$ is bounded on bounded sets; on the other hand, the difference
$$
(-\Delta)^s\phi(F_1)-(-\Delta)^s\phi(F_2)=O(|x|^{-(1+2s}) \quad \mbox{as } \ |x|\to\infty.
$$
The last formula comes directly from the representation formula for $(-\Delta)^s$ plus the fact that $\phi(F_1)-\phi(F_2)$ has compact support.  This correction  $O(|x|^{-(1+2s)})$ is lower order with respect to $(-\Delta)^s\phi(F_1)$ at infinity since $(-\Delta)^s\phi(F_1)\sim F_1=F_2\sim |x|^{-2s/(1+n)}$ for all large $|x|$.  Summing up, there exist constants $K_1, K_2>0$ such that
\begin{equation}
-K_1 F_2(x)\le -(-\Delta)^s\phi(F_2)\le K_2 F_2(x)\,.
\end{equation}
If we now take  $H(0)=a>0$ and $H'\le -K_1 H^{-n}$, we  get a subsolution,
 $$
 \partial_t \widetilde U+(-\Delta)^s \Phi_n(\widetilde U)\le 0\,,
 $$
in some time interval. In the last case maybe $H$ decays and vanishes in finite time (even if this is not realistic for actual solutions, as we will see, it is just the form of the modified subsolution). Take $F=F_2$ and this choice of $H$ to end the construction of the desired subsolution $\widetilde U$.  Finally, note that since $H(0)$ can be taken at will we find a family of bounded subsolutions and the $L^\infty$ norm can be taken as small as wanted.

\medskip

(ii) In the same way, if we take $F_1$ as before and $H(0)=a>0$ such that $H'\ge K_2 H^{-n}$, then we get a formal supersolution. In this case it will exist for all times. \qed

\medskip

\noindent $\bullet$ {\bf Case $n=0$.}  This case is settled by replacing the power $-u^{-n}$ by $\log(u)$ and repeating the above procedure. We will need the calculation the $s$-Laplacian of the logarithm. We have
\begin{equation}\label{slap.log}
(-\Delta)^s (\log|x|)=c(s) \,|x|^{-2s}\,.
\end{equation}
A short proof is as follows: using the previous formula for $\alpha>0$, $\alpha$ very small we get
\begin{equation}
k(\alpha,s)\sim  \alpha 2^{2s-2}\frac{(2s-1)\Gamma(1/2)\,\Gamma(s)}{\Gamma((3-2s)/2)}
\end{equation}
where we have used $\Gamma(-\alpha/2))\sim \Gamma(1)(-2/\alpha)$, and $\Gamma((1+\alpha-2s)/2)\sim 2/(1-2s)\,\Gamma((3-2s)/2)$. We now use the expression $\log(x)=\lim_{\alpha\to 0} (x^\alpha-1)/\alpha$, $x>0$, to conclude that the formula is true with $c(s)=\lim_{\alpha\to 0} k(\alpha,s)/\alpha>0$ if $2s>1$.

Then the rest of the steps is quite similar and the conclusions of Lemma \ref{lemma.subsol}  hold.
\qed

\medskip

\begin{prop} \label{prop.subsol} Let  $N=1$ and $0< n< 2s-1$. Let the initial data $u_0$ be positive and integrable and satisfy $u_0(x)\ge C F_2(x)$ for all $x$. Then the limit solution constructed as in Section {\rm \ref{sec.limit}} is non-trivial. In fact, it sits on top of one of the constructed subsolutions, hence it is a positive very weak solution of the equation, at least in a certain time interval \ $0<t<T$.
Nontrivial solutions are also obtained for \ $n=0$, $1/2<s<1$.
\end{prop}

\noindent {\sl Proof.}  The subsolution can be compared by classical results with all approximate problems and remains below for all the existence time. This allows us to prove that the limit solution  does not vanish identically when we take initial data that decay equal or slower than $F$. By the $L^1$ contraction the same nontrivial limit happens for any continuous and bounded initial data.

\medskip

\noindent {\bf Remark.} This is the first successful step in a long road that leads to the proof that all upper limits corresponding to nontrivial data are in fact positive weak solutions.


\subsection{Some properties of nontrivial limit solutions}\label{limit.proprts}

\noindent  $\bullet$ {\bf Scaling property.} If $u(x,t)$ is a nontrivial limit solution, then so is
\begin{equation}\label{scale}
 u_{AL}(x,t)= A u(Lx, L^{2s}A^{-(1+n)} t)
\end{equation}
for all  parameters $A,B>0$. We leave the easy proof to the reader (see similar arguments in \cite{JLVSmoothing}).

\medskip

\noindent  $\bullet$ {\bf The smoothing effect}. It says that all limit solutions with integrable data are in fact bounded so that we use bounded solutions in the proofs. We have for all limit solutions
\begin{equation}\label{smoothing.L1}
\sup_{x\in{\mathbb{R}}}|u(x,t)|\le
C_2 \,t^{-\alpha}\|u_0\|_{1}^{\delta}
\end{equation}
with
$\alpha=1/(2s-n-1)$ \ and \ $\delta=2s/(2s-n-1)$,
the constant $C_2$ depending on $n,s$.

\noindent {\sl Proof.} Use the smoothing effect in the rough form already proved for the approximate solutions and let $t=1$ to conclude that the result holds for $M=\|u_0\|_{1}=1$. When $M\ne 1$ and $t\ne 1$ use the scaling rule \eqref{scale} in the usual way to reduce the proof to the particular case. This is a well-known scaling trick. \qed

\noindent  $\bullet$ {\bf Time monotonicity.} The limit solutions satisfy
\begin{equation}\label{monot.time}
\partial_t u\le \frac{u}{(1+n)t}\,.
\end{equation}
It follows from the same property for the approximations. As we said there, the argument has a proof using scaling arguments originally due to B\'enilan and Crandall \cite{BenCr}.

Also the properties of $L^p$ boundedness, pointwise comparison, and $L^1$ contraction of the approximate solutions pass to the limit without change.

\noindent  $\bullet$ {\bf Space monotonicity. Aleksandrov's principle.}
The Aleksandrov-Serrin reflection method is a well-established tool to prove monotonicity of solutions of wide classes of (possibly nonlinear) elliptic and parabolic equations,  cf. \cite{Alek60, Ser71}. It has been quite useful in particular in the case of the PME, as documented in  \cite{CVW87, Vapme}. This is the version proved in \cite[Theorem 15.2]{VazBar} for nonlinear parabolic equations with fractional diffusion of the type \eqref{geq.Phi}, and adapted to our situation

\begin{prop} Let $v_\eps$ the unique solution of \eqref{eq:FrUF_approx1} with initial data $u_0\in L^1(\re)$, $u_0\ge 0$. Under the assumption that
\begin{equation}
u_0(x)\le u_0(2a-x)) \qquad \text{for} \quad x>a
\end{equation}
for some $a>0$, we have for all $t>0$
\begin{equation}
v_\eps(x,t)\le v_\eps(2a-x,t) \qquad \text{for} \quad x>a\,.
\end{equation}
\end{prop}
In plain words, the result deals with comparison of a solution with its space reflection with respect to the point $x=a$. If it is true for $t=0$, then it is true forever. After passage to the limit the same comparison is true for the obtained upper solutions. An immediate consequence of this
that has been used in the literature and we will use below is

\begin{cor}\label{aleks} An upper limit solution with initial data supported in the half-fine $\{x<a\}$ is monotone nonincreasing in $x$ in the region $\{x>a, t>0\}$. If the initial data are supported in the half-fine $\{x>-a\}$ the solution is monotone nondecreasing in $x$ in the region $\{x<-a, t>0\}$.
\end{cor}


\section{Comparison results}\label{sec.comp}

The standard comparison theorem (Maximum Principle) applies the approximate problems, hence it will be valid in the limit $\ve\to 0$ for the class of upper limit solutions. The Aleksandrov principle is another comparison theorem. In the sequel we will use two other comparison results, that we discuss next.

\subsection{Symmetrization and concentration comparison}\label{ssec.sym}

 Symmetrization techniques are a very popular tool of obtaining a priori estimates for the solutions of different partial differential equations, notably those of elliptic and parabolic type. The application of Schwarz symmetrization allows to obtain sharp a priori estimates for elliptic problem by comparison with a model symmetric problem. For parabolic problems the usual pointwise comparison of the solutions of the two problems fails, and is replaced by comparison of integrals, \cite{Bandle}. In the case of the porous medium equation $u_t=\Delta u^m$ that result was established in \cite{Vsym82, VANS05}, and holds for all $m>0$. In order to state the result we will use, the following definition is needed:

\noindent {\bf Definition.} Let $f,g\in L^{1}_{loc}(\ren)$ be two radially symmetric functions on $\ren$. We say that $f$ is less concentrated than $g$, and we write $f\prec g$, if for all $R>0$ we get
\begin{equation}\label{form.cc}
\int_{B_{R}(0)}f(x)\,dx\leq \int_{B_{R}(0)}g(x)\,dx.
\end{equation}

The partial order relationship $\prec$ is called \emph{comparison of mass concentrations}.
In the applications we are going to assume that $f$ and $g$ are  {\sl rearranged} functions.

The following result is proved in \cite{VazVol1}.

\begin{theor} Let Let $u_1, u_2$ be two nonnegative, weak solutions of  the equation  $u_t+ (-\Delta)^s\Phi(u)=0,$ posed in $Q=\ren\times(0,\infty)$, with nonnegative initial data $u_{01}, u_{02}\in L^1(\ren)$. Assume that both $u_{02}$ and $u_{01}$ are rearranged and  $u_{02}\prec u_{01}$. Assume moreover that the nonlinear function $\Phi(u)$ is positive, smooth and concave for $u>0$. Then, for all $t>0$ the functions $u_1(\cdot,t)$ and $u_1(\cdot,t)$ are rearranged and we have
\begin{equation}\label{form.cc2}
u_2(\cdot,t)\prec u_1(\cdot,t).
\end{equation}
In particular, we have $\|u_2(\cdot,t)\|_p \le\|u_1(\cdot,t)\|_p$ for every $t>0$ and every $p\in [1,\infty]$.
\end{theor}

This result applies to the solutions $v_\ve$ of the regularized problems \eqref{eq:FrUF_approx}. In the limit it will apply to all nontrivial solutions of the singular fractional FDE \eqref{geq.Phi} with $\Phi=\Phi_n$ and $n\ge 0$, for all $0<s<1$.

Note that in 1D the integrals in formula \eqref{form.cc} have a more classical interpretation: the integral $\int_0^x f(x)\,dx$ is just the distribution function of the mass density $f$, that we are assuming to be nonnegative, integrable and monotone decreasing for $x>0$. Therefore, formula \eqref{form.cc2} is just a comparison of distribution functions.

\medskip

\subsection{Shifting comparison}

A new comparison result that is related to symmetrization in spirit and techniques is based on lateral displacement of the solution, viewed as a mass distribution.

\begin{theor}\label{thm.shift} Let two functions $u_{01}, u_{02}\in L^1(\re)$ with the following properties:

(i) They are nonnegative and rearranged around their points of maximum $x_1$ and $x_2$ resp., with $x_1<x_2$. The total mass is the same.

(ii) We assume moreover that $\int_{-\infty}^x u_{02}\,dx\le \int_{-\infty}^x u_{01}\,dx$ for every $x\in\re$.

\noindent Besides, we assume that $\Phi$ is monotone, concave and defined on $\re_+$ with $\Phi(0)=0$ and $0<\Phi'(0)<\infty$. Then,  the following comparison inequalities hold for the corresponding (limit) solutions
\begin{equation}
\int_{-\infty}^x u_2(x,t)\,dx\le \int_{-\infty}^x u_1(x,t)\,dx,
\end{equation}
for every $x\in\re$ and every $t>0$.
\end{theor}

This result is called the {\sl Shifting comparison} lemma. It is essentially one-dimensional and it was established in the PME case by V\'azquez \cite{VTr1} and  it proved useful in studies of free boundary location or asymptotic behaviour.  It is related to mass transport and Wasserstein distances, \cite{Villani2003, VazGaeta}. It will be crucial in some proofs below, like the proof of the existence for general initial data and the asymptotic behaviour. Since it has  a rather technical and long proof, we will delay to Section \ref{sec.shifting} at the end of the paper.


\section{\bf Mass conservation and global solutions}\label{sec.cm1}

The property of mass conservation plays an important role in passing from local-in-time existence
to global solutions, since it prevents the phenomenon of finite-time extinction. A first result is as follows.

\begin{prop}\label{prop4.1} Let us assume that $N=1$, $s>1/2$ and \ $0\le n<2s-1$ and let us assume that $u_0$ is integrable and $u_0(x)\ge cF_2(x)$, i.e.,
\begin{equation}\label{bdd.below}
u_0(x)\ge c/(1+\,|x|^2)^{s/(1+n)}
\end{equation}
for some $c>0$. Then the limit solution not only is positive in a certain time interval $0<t<T$, but it also conserves mass in that interval:
\begin{equation}
\int_{\re} u(x,t)\,dx=\int_{\re} u_0(x,t)\,dx.
\end{equation}
\end{prop}

\noindent {\sl Proof.} (i) Let us first assume that $n>0$. We will prove conservation of mass for small times. We take a nonnegative non-increasing cut-off function $\zeta(s)$ such that $\zeta(s)=1$ for $0\le |s|\le1$, $\zeta(s)=0$ for $|s|\ge2$, and define
$\zeta_R(x)=\zeta(|x|/R)$.  We have $(-\Delta)^{s}\zeta_1\in L^1(\mathbb{R}^N)\cap
L^\infty(\mathbb{R}^N)$, moreover, $|(-\Delta)^{s}\zeta_1|\sim |x|^{-(1+2s)}$ as $|x|\to \infty$. The radial cut-off function $\zeta_R$ has the scaling property
\begin{equation}
\label{varphi1}
(-\Delta)^{s}\zeta_R(x)=R^{-2s}(-\Delta)^{s}\zeta_1(x/R).
\end{equation}
We take the approximate solutions $u_\ve$ of equation \eqref{geq.Phi} with nonlinearity $\Phi_\ve$ as in Section \ref{sec.limit}, multiply by  $\zeta_R$ and integrate by parts. We have
$$
\int_{\mathbb{R}} (u_\ve(t_1)-u_0-\ve)\,\zeta_R\,dxdt=
\int_0^{t_1}\int_{{\mathbb{R}}}u_\ve^{-n}\,(-\Delta)^{s}\zeta_R\,dx.
$$
Note that the last integral is absolutely integrable.
Passing to the limit $\ve \to 0$ we get for every $t>0$,
\begin{equation}\label{id.cm1}
\int_{\mathbb{R}} u(t_1)\,\zeta_R\,dx -\int_{\mathbb{R}} u_0\,\zeta_R\,dx=
\int_0^{t_1}\int_{{\mathbb{R}}}u^{-n}\,(-\Delta)^{s}\zeta_R\,dxdt.
\end{equation}
Let us split the right-hand side of \eqref{id.cm1}
into the integrals for $|x|\le R$ and $|x|\ge R$. We only estimate the integrals in $x$, forgetting for the moment the time integration. We get
$$
|I_1(R)|\le \int_{|x|\ge R} u^{-n}(t)\,|(-\Delta)^{s}\zeta_R||\,dx
\le \frac{C}{R^{2s}} \int_{|x|\ge R} |x|^{2ns/(1+n)}(|x|/R)^{-(1+2s)}\,dx.
$$
Putting $x=Ry$ we get
$$
I_1(R)\le \frac{CR^{2ns/(1+n)}}{R^{2s-1}}\int_{|y|\ge 1} |y|^{\frac{2ns}{1+n}-(1+2s)}\,dy=
CR^{-\gamma}\int_{|y|\ge 1}\frac{dy}{|y|^{2+\gamma}}
$$
Since $\gamma= 2s-1-2ns/(1+n)=2s/(n+1)-1>0$, we get $I_1(R)\to 0$ as $R\to\infty$. On the other hand,
for the analogous integral in the set $|x|\le R$ we get
$$
|I_2(R)|\le \int_{|x|\le R} u^{-n}(t)\, |(-\Delta)^{s}\zeta_R|\,dx\le
 \frac{C}{R^{2s}} \int_{|x|\le R} |x|^{2ns/(1+n)}\,dx\le CR^{-\gamma}
$$
that goes also to zero as $R\to\infty$. Since these estimates do not depend on $t$ (for small $t$)
we may go back to equation \eqref{id.cm1} and let $R\to\infty$ to get
$$
\int_{\mathbb{R}} u(x,t_1)\,dx =\int_{\mathbb{R}} u_0(x)\,dx,
$$
which is the mass conservation law. This law holds for the small times for which we have the lower estimate for the limit solution $u$ (for instance, the estimate coming from comparison with the subsolution constructed in Section \ref{sec.bddsub}, as used in Proposition \ref{prop.subsol}).

This method of proof is inspired in a common technique that was used to prove the property in the  case $s=1$ for $m>(N-2)/N$, $N\ge 2$. The beginning is the same, but the further details are quite different.

\medskip (ii) In the case $n=0$ the nonlinearity is logarithmic and not power-like. We arrive at
formulas $I_1(R)$ and $I_2(R)$ with the bound $c\log(1+|x|^2)$ instead of $c|x|^{2ns/(1+n)}$, and we put $n=0$ in the rest of the places. The remaining steps follow easily.  \qed

 \medskip

\begin{prop}\label{prop4.2} If the initial data $u_0$ is a rearranged function with mass $M>0$, the limit solution $u(x,t)$  exists globally in time, is positive everywhere and the mass is conserved for all times.
\end{prop}

\noindent {\sl Proof.} (i) We first prove that for  rearranged initial data in the same class $u_0\ge c\,F_2(x)$ the solution exists for all times and conservation of mass holds also for all times. We know  that $u$ is positive and conserves mass for $0<t<T_1=T_1(u)$. Assume that this maximal time is finite. We use  the scaling property of the equation to define of new solutions
\begin{equation}\label{mc.scaling}
 u_L(x,t)= L u(Lx, L^{2s-(1+n)} t)
\end{equation}
with $L>1$.  Let us introduce the notation $u_L={\cal T}_L u$ for future reference. This scaling, a particular case of \eqref{scale}, keeps the mass of the solutions $u$ and $u_L$ identical. Now, since the time of existence of the family $u_L$ shrinks to the time $T_L=T_1/L^{2s-(1+n)}$, we seem have a problem  in using rescaling. But the problem can be fixed by using  symmetrization (concentration comparison) and we end up with an expected gain.

Indeed, the new solution $u_L$ for $L>1$ is clearly more concentrated than $u_1(x,t)=u(x,t)$ at at time $t=0$, hence it will be more concentrated at all times  by the result of Subsection \ref{ssec.sym}. Then, the concentration relation \eqref{form.cc2} immediately implies that the mass of $u_L$ must be conserved as long as the mass of $u$ is, say until $T_1(u)$; the justification for $T_L\le t\le T_1$ is done by doing symmetrization comparison on the approximate problems and passing to the limit, that cannot be trivial because of this argument. This means that the nontrivial existence time  with conservation of mass for  $u_L$ is $T(u_L)\ge T_1$. Undoing the scaling  we get the same property for $u$ in a time $T(u) \ge L^{2s-(1+n)}T_1$, hence $T_1$ must be infinite.

\medskip

\noindent (ii) Next, we prove that the solution is positive everywhere.
Let us do an analysis of what happens if a rearranged  solution touches zero, and show that this cannot happen. First, we use the monotonicity properties in space and time to show that, if the solution vanishes at $t=t_1$ and $x=R$, then we must have $u(x,t)=0$ for all $|x|\ge R$ and $t\ge t_1$. Using the proof  of the mass formula of Proposition \ref{prop4.1} we have for the approximations $u_\ve$
$$
\frac{d}{dt}\int_\re u_\ve(x,t)\zeta(x)\,dx=-\int_\re \Phi_n(u_\ve)\,(-\Delta)^s \zeta\,dx
$$
Taking a cutoff function supported in  $[-R_1,R_1]$, with $R_1<R$, we know that $-(-\Delta)^s \zeta>0$ for $|x| \ge R_1$, and since $\Phi_n(u_\ve)$ tends to minus infinity in the set we have described, we have
$$
-\int_{|x|\ge R_1} \Phi_n(u_\ve)\,(-\Delta)^s \zeta\,dx\to -\infty
$$
as $\ve\to 0$. The integral for $|x|\le R_1$ has a bounded limit since the limit $u$ is bounded below, hence $|\Phi(u)|$ is bounded. Applying this between times $t_1$ and $t_1+\tau$, we can show that the weighted mass must be zero for all times larger than $t_1$. Since the solution is rearranged, this implies that  the mass will be zero, which is excluded by the previous step. \qed

\noindent (iii) The next step is to consider rearranged initial data that do not satisfy a bound from below like \eqref{bdd.below}, for instance $u_0$ may be compactly supported. Let $M>0$ be the initial mass. In that case we make a small perturbation by adding to $u_0$ a tail of the form $\eta^\delta(x)=\delta (1+ |x|^2)^{-s/(1+n)}$ and besides we truncated the initial data on top to make it bounded. We may do all this and conserve the mass $M$. In this way we obtain a solution $u^\delta$ to the problem where the previous analysis applies, and mass is conserved. We now use the $L^1$ contraction property and the conservation of mass for $u^\delta$, we conclude that the limit solution $u$ corresponding to initial data $u_0$ must have mass
$$
\int u(x,t)\,dx\ge M-\|u_0^\delta-u_0\|_1
$$
which is positive for $\delta$ small, hence $u$ is a global solution. The argument is justified as limit solution, i.\,e., in the  limit of approximate problems. Finally, by letting $\delta\to 0$ we derive the mass conservation property for $u$. \qed

\medskip

Positivity can now be obtained for all solutions with continuous initial data, not necessarily rearranged, by standard comparison with a solution with compactly supported rearranged data (after possibly a space displacement to fit it under $u_0$). Such a subsolution  is positive for all positive times, hence $u$ is too. We will return to this question later on, after we get some quantitative estimates on the behaviour.

\subsection{Concept of  solution. The mass function}
As a preliminary  for the next developments, we need to clarify the type of solution that we get at this stage when we pass to the ``limit solutions''. In the end, we would like to prove that our positive solutions are very weak in stated sense that
\begin{equation}
\int\int u\,\zeta_t\, dxdt =\int\int \Phi_n(u)\,(-\Delta)^s \zeta\,dxdt
\end{equation}
for all smooth test functions with compact support. However, the last term offers a difficulty as long as we do not know the decay of $u$ at infinity, i.\,e., as long as we do not control the growth of $u^{-n}$, since typically $(-\Delta)^s \zeta$ behaves like $O(|x|^{1+2s})$ as $|x|\to\infty$.

The conditions of Proposition \ref{prop4.1} do allow for a correct passage to the limit $\ve\to0$ in the definition, but the conditions of Proposition \ref{prop4.2} do not. Keeping the assumption of rearranged data, we find a  remedy by weakening the definition by integration in space to get a new function $V(r,t)=\int_0^r u(x,t)\,dx$ so we have $V_r=u$ for all $r>0$ and we can write the equation formally as
\begin{equation}
V_t=\int_0^r u_t\,dx= -\int_0^r (-\Delta)^s \Phi_n(u)\,dx
\end{equation}
Putting $(-\Delta)^s= -\partial_{xx}^2(-\Delta)^{s'}$ with $s'=1-s>0$, and integrating we get
 \begin{equation}\label{eq.mass}
V_t= \partial_r(-\Delta)^{-s'} (\Phi_n(u)),    \qquad V_r=u\,.
\end{equation}
This is an integrated version of the weak solution that makes perfect sense for the approximate problems and for their limits in the very weak sense. We will call $V(x,t)$ the {\sl mass function; } it is the  distribution function in Probability, but that name might lead to confusion here.


\section{The Barenblatt solutions in 1D}\label{sec.bar1d}

After establishing mass conservation for all times we can construct the Barenblatt solutions
and derive the main properties.

\subsection{Existence}
 Take one of the rearranged solutions $u_1(x,t)$ of the previous section (Proposition \ref{prop4.2}) with initial data $u_{01}(x)$ that we may assume continuous. It exists globally in time, is positive everywhere and conserves mass. Let us fix the $L^1$ norm of $u_1$ to 1. Take then the rescaled family $\{u_L(x,t)={\cal T}_L u:\ L\ge 1\}$ defined by formula  \eqref{mc.scaling} of the previous section. Finally,  pass to the limit
\begin{equation}
U(x,t)=\lim_{L\to \infty} u_L(x,t).
\end{equation}
We have to show that this limit exists and has the desired properties.

(i)  In principle, the family $u_L(x,t)$ converges weakly in $L^1$ for $t\ge \tau>0$, and there may also be a non-unique limit. It is best to use an argument based on concentrations, that makes it natural to argue with the family of mass functions $V_L(x,t)=\int_0^x u_L(x,t)\,dx$. By the concentration comparison result we see that the family $\{V_L\}$ is  monotone increasing in the parameter $L$ for $x>0$ (resp. negative and decreasing for $x<0$). Hence, the limit $V_\infty(x,t)$ exists and does not depend on subsequences $L_k\to\infty$. The convergence $V_L(x,t)\to V_\infty(x,t)$ is uniform convergence in $x$ for every fixed positive time.

(ii) Differentiation in $x$ gives the unique weak limit $U=\lim_{L\to\infty} U_L$. Stronger convergence will be proved later on.

(iii) For $t=0$ we have $V_\infty(x,0)=1/2$ for all $x>0$, $V_\infty(x,0)=-1/2$ for all $x<0$, in other words the limit of the initial data is a delta function. We have to show that for all positive times $U(\cdot,t)$  is not trivial  away from $x=0$ (i.e., it is not equal to $\delta(x)$).  This follows from the smoothing effect (proved in Subsection \ref{limit.proprts}) that applies uniformly to all functions $u_L$. We conclude that $U(\cdot,t)$ is a bounded function for all $t>0$, therefore $U(\cdot,t)$ is not a Dirac delta, and $V_\infty(\cdot,t)$ is  Lipschitz continuous in $x$ uniformly if $t  \ge \tau>0$.

(iv) Let us prove next that the limit must be a self-similar function. The argument is based on passing to the limit in the scaling family using the group of transformations. In fact, for all $L,k>0$ we have ${\cal T}_L {\cal T}_k u_1={\cal T}_{kL}u_1$,
Passing to the limit $k\to\infty$ and using the the uniqueness of the limit because of the comparison of concentrations, we get $U(x,t)={\cal T}_L U(x,t)=L U(Lx, L^{2s-1-n}t)$, hence $L U(Lx, L^{2s-1-n})=U(x,1)$. In the usual way it follows that
\begin{equation}\label{f.bar.1}
U(x,t)=t^{-\alpha}F(xt^{-\alpha}), \quad \alpha=1/(2s-(1+n)\,,
\end{equation}
with $F(x)=U(x,1)$.

(v) It is immediate that the profile $F$ is positive everywhere, bounded, integrable, and rearranged (i.\,e., $F(x)=F(-x)$ and $F$ monotone decreasing for $x>0$).

(vi) To construct the Barenblatt solutions with initial data $M\,\delta(x)$ with any mass $M>0$ we use another scaling ${\cal T}'$ defined by $({\cal T}'u)(x,t)=M\,u(x,M^{-(1+n)}t)$, that transforms solutions of mass 1 into solutions of mass $M$. We get a self-similar solution with the same formula as before, but now the profile is
\begin{equation}
F_M(x)=M^{2s\alpha}F(xM^{(1+n)\alpha})\,.
\end{equation}

\subsection{Alternative approach and better convergence}

The idea is to apply the same type of approach to the approximate problems \eqref{eq:FrUF_approx} to obtain a fundamental solution of each of those problems. After careful inspection, we see that we can find a solution $U_{\infty\eps}(x,t)$ with initial data $U_{0\eps}(x,0)=\delta_0(x)+\eps$. This is better done by using the formulation $v_\eps(x,t)=u_\eps(x,t)-\eps$ that solves Problem \eqref{eq:FrUF_approx1}, to which we can apply the usual $L^1$ theory and
comparison of concentrations. Therefore, we obtain a fundamental solution with this argument, the difference is that this time we cannot conclude that it is self-similar. On the other hand, the whole collection of rescaled solutions $U_{L,\eps}(x,t)$ are uniformly bounded for $t\ge \tau>0$ hence, by the regularity results of \cite{VPQR13}, they are uniformly $C^\alpha $ continuous for a certain $\alpha>0$. The family is thus locally compact in $L^1_{loc}$ (in both space and time), which means that the convergence $U_{L,\eps}(x,t)\to U_{\infty,\eps}(x,t)$ takes place in the strong sense of $L^1(\re)$ for every $t>0$. This is a key improvement in the situation.

Now we take the monotone limit of these fundamental solutions to get
$$
\lim_{\eps\to 0}U_{\infty,\eps}(x,t)=U_\infty(x,t)
$$
By a simple comparison, $U_\infty(x,t)\ge U(x,t)$, where $U$ is the previously constructed Barenblatt solution, formula \eqref{f.bar.1}. By the equality of masses we conclude that both are the same function, $U_\infty=U$.

It is now easy to see that $U_L(x,t)\to U(x,t)$ in $L^1(\re) $ for every $t\ge \tau>0$. In fact,
\begin{align*}
\|(U_L(x,t)-U(x,t))_+\|_{L^1(-R,R)}\le \|(U_{L,\eps}(x,t)-U(x,t))_+\|_{L^1(-R,R)}\\
\le \|(U_{L,\eps}(x,t)-U_{\infty,\eps}(x,t))_+\|_{L^1(-R,R)}+ o(\eps),
\end{align*}
where the last term is the contribution of $\|U_{\infty,\eps}(x,t)-U_{\infty}(x,t)\|_{L^1(-R,R)}$. Fixing a small $\eps>0$ and and using the convergence $U_{L,\eps}(x,t)\to U_{\infty,\eps}(x,t)$ in $L^1(\re)$, there is an $L_\eps>0$ such that have $\|U_L(x,t)-U(x,t))\|\le \eps$ for $L\ge L_\eps$. All together, this means that
 $(U_L(x,t)-U(x,t))_+\to 0$ in $L^1(-R,R)$ as $L\to\infty$.  Since the mass in the far field is uniformly small, and there is total mass equality, we get $U_L(x,t)-U(x,t)\to 0$ in $L^1(\re)$. This is the improved convergence that we needed.

\subsection{Uniqueness}

\begin{prop} The fundamental solution is independent of the rearranged function $u_0$ that starts the construction.
\end{prop}

\noindent {\sl Proof.} For bounded functions of compact support this is done by symmetrization comparison and iterated limits. Indeed, let us start from functions $u_1$ and $\tilde u_1$ and let $U$ and $\tilde U$ the corresponding Barenblatt solutions obtained in the limit of the scaling process. Then it is not difficult to see that the scaled function $u_{0k} $ is more concentrated than $\tilde u_{0}$ if $k$ is large enough. Therefore, $u_{kL}$ is more concentrated than $\tilde u_{0L}$, hence in the limit $U(\cdot,t)$ is more concentrated than $\tilde U(\cdot,t)$. The reverse relation also holds, hence $U=\tilde U$.

For general data we use $L^1$ contraction. \qed

\section{Study of the profile}\label{sec.F}

We want to know more about the profile $F$ of the Barenblatt solution, in particular its equation and its asymptotics for large $|x|$.

\subsection{Profile Equation}
We  can apply to $U$ the concept of solution of last section and perform the computation when $U$ is selfsimilar. Then, with $\xi=x\,t^{-\alpha}$,
$$
V_t=\frac{d}{dt}\int_0^r t^{-\alpha}F(x\,t^{-\alpha})\,dx=-
\alpha\,t^{-(\alpha+1)}\int_0^r (F(\xi)+\xi\,F'(\xi))\,dx=
-\alpha\,t^{-1}\xi\,F(\xi)\,.
$$
Therefore, equation \eqref{eq.mass} becomes at $t=1$
\begin{equation}\label{eq6.3}
\partial_x(-\Delta)^{-s'} (F^{-n}(x))=\alpha x \,F(x), \quad x>0.
\end{equation}
This is the  integro-differential equation satisfied by the profile $F$.  Note that $1-2s'=2s-1>0$, so that the operator in the LHS has a positive degree of differentiation, $\partial_x{\cal L}^{-s'}=H{\cal L}^{s-1/2}$. We can also write the equation as
\begin{equation}\label{form.7.2}
\int_0^r (-\Delta)^s F^{-n}(r)\,dr=\alpha  x \,F(x), \quad  \int_0^r  F^{-n}(r)\,dr=\alpha (-\Delta)^{-s }(x \,F(x))\,.
\end{equation}

\subsection{Asymptotic behaviour of $F$}
According to the analogy with the Barenblatt solutions constructed for the case $m>0$ in  \cite{VazBar}, since we are in the case $m<m_1=N/(N+2s)$ (a value that was important in that respect), we expect the following behaviour.

\begin{prop}\label{precise.ab} The profile $F$ decays as $|x|\to\infty$ like $r^{-\gamma}$, with $\gamma=  2s/(1+n)>1\,$. Indeed, there exists the finite positive limit
\begin{equation}
\lim_{r\to+\infty} F(r)\,r^{\gamma}=c_\infty\,,
\end{equation}
and moreover $F(r)\le c_\infty r^{-\gamma}$  for all $r>0$.
\end{prop}

\noindent {\bf Remark.} Later on we will calculate   the constant $c_\infty$ explicitly as a corollary of our work on very singular solutions, Section \ref{sec.vss}.

\medskip

\noindent {\sl Proof.} (i) We need a first estimate of the decay of $F$. Since it is an integrable and rearranged function we immediately get $F(x)\le C/|x|$, which is too rough.  A better estimate is obtained as follows: We start the scaling procedure to construct $U$ by taking as initial data $u_{01}$ the function $F_2(x)$ of Lemma \ref{lemma.subsol}. By comparison with the subsolution of the Lemma we have a decay rate $u_1(x)\ge c(t)|x|^{-2s/(1+n)}$ for all $|x|\ge 1$ uniformly in some time interval $0<t<T_1$.

We can show a similar decay rate for the limit solution $U(x,t)$ by applying concentration comparison. Indeed, the mass of $u_1$ in region $|x|\ge R>1$ is estimated as
$C(t)R^{1-2s/(1+n)}$, and the mass of $U$ has to be less than that by the comparison. By using the monotonicity of $U$ w.e.t. $|x|$ we get the conclusion that $U$ is less than $c_1(t)\,|x|^{-2s/(1+n)}$. Recall that $F(x)$ is just $U(x,1)$.  Therefore,
$$
F(x)\le C\,|x|^{-2s/(1+n)}\,.
$$
The reader can find similar arguments in \cite[Section 12]{VazBar}.

(ii) The power-like bound from below can be obtained following the ideas of \cite{VazBar}: We start from the homogeneity estimate \eqref{monot.time} that says that $(1+n)\,t\,u_t\le u$. In terms of the self-similar profile,  this just means that $-(1+n)\alpha(\,F + rF'(r))\le F$, $r=|x|>0$, hence
$$
\frac{-rF'(r)}{F(r)}\le 1+\frac{1}{(1+n)\alpha}=\frac{2s}{1+n}\,.
$$
Integration of this inequality gives the following lower bound, valid for all $r\ge 1$, all $s\in (0,1)$:
\begin{equation}\label{fde.rate}
F(r)\ge C\, r^{-2s/(1+n)}\,.
\end{equation}
Moreover, the function $ J(r):= F(r)\, r^{-2s/(1+n)}$ is monotone non-decreasing with $r$, so that it has a limit as $r\to\infty$. Let us put
\begin{equation}
\lim_{|x|\to \infty } F(x)\,|x|^{-2s/(1+n)}=c_\infty(s,n)\,.
\end{equation}
In principle the limit may be finite and positive or infinite. By part (i) it is finite.

\section{The very singular solution}   \label{sec.vss}

The a priori bounds on the profile of the Barenblatt solutions make it easy to pass to the limit $M\to \infty$ and obtain a special function, called the very singular solution (VSS), much as we have done in \cite{VazBar} for $m>0$. This is the result

\medskip
\begin{theor}\label{thm.vss2}
 The Very Singular Solution calculated in Section \ref{sec.bddsub},
 \begin{equation}
 {\widetilde U}(x,t)= C(n,s)\,t^{1/(1+n)}|x|^{-2s/(1+n)},
 \end{equation}
   is the limit of the Barenblatt solutions $U_M(x,t)$ as the mass $M$ goes to infinity.
\end{theor}

\noindent {\sl Proof.} We recall that the value of the constant $C(n,s)$  can be explicitly computed from the  calculations there  as $C(n,s)=K(n,s)(1+n)$.

By the established comparison properties, it is clear that the sequence of Barenblatt solutions $U_M(x,t)$ is monotone increasing with $M>0$. Next, we check that they are all bounded above by the VSS. A direct comparison of $U_M$ with ${\widetilde U}$ is difficult to justify directly, hence we argue in another way to get our conclusion.

We first fix the mass equal to one and use the upper profile bound, $F(r)\le c_\infty r^{-2s/(1+n)}$  for all $r>0$, to conclude that
$$
U_1(x,t)\le c_\infty t^{1/(1+n)}|x|^{-2s/(1+n)}\,.
$$
We know that the whole sequence $U_M$ can be obtained from $M=1$ by the rescaling
$U_M(x,t)= ({\cal T}_M u)(x,t):=M\,u(x,M^{-(1+n)}t)$. We immediately see that $U_M(x,t)$ satisfies the same upper bound, even with the same constant, $c_\infty(M)=c_\infty(1)$. Once, we have the same upper bound for the whole sequence, we may pass to the monotone limit and get
$$
U_\infty(x,t)\le c_\infty t^{1/(1+n)}|x|^{-2s/(1+n)}\,.
$$
But since the functions $U_M$ are invariant under the mass conserving scaling,  so is the unique limit, hence
$U_\infty(x,t)$ is self-similar, $U_\infty(x,t)=t^{-\alpha}F_\infty(x\,t^{-\alpha})$. Also $F_\infty(y)\ge F_M(y)$ for all $y>0$ and all $M>0$, and $F_\infty(y)\le c_\infty |y|^{-2s/(1+n)}$. It easily follows from a tail analysis that $F_\infty(y)=c_\infty |y|^{-2s/(1+n)}$, hence
$$
U_\infty(x,t)= c_\infty t^{1/(1+n)} |x|^{-2s/(1+n)}\,.
$$
We conclude that $U_\infty$ is another possible very singular solution of the equation obtained as limit of upper limit solutions. In order to see that $U_\infty$ and $\widetilde U$ must be the same we only have to check that both are weak solutions for the equation for $x \ne 0$, in other words, the profiles must be weak solutions for the $F$ equation \eqref{eq6.3}. This is what selects the constant in a unique way. We conclude that $c_\infty=C(n,s)$ and $U_\infty=\widetilde U$. \qquad

As a corollary of the above result we can refine the information on the Barenblatt solutions obtained in Proposition \ref{precise.ab} as follows
\begin{prop} \label{precise.ab.2} The profiles $F_M$, $M>0$, of the Barenblatt solution satisfy
the uniform bound $F(r)\le C(n,s) r^{-2s/(1+n)}$  for all $r>0$ and this estimate is sharp at infinity
\begin{equation}
\lim_{r\to+\infty} F(r)\,r^{2s/(1+n)}=C(n,s)    \,,
\end{equation}
where $C(n,s)$ is the constant of the VSS.
\end{prop}


\section{The theory for general initial data}

Here we want to complete the proof of Theorems \ref{thm.1} and \ref{thm.2} on the existence and properties of solutions of the Cauchy Problem. We establish existence, positivity and behaviour as $|x|\to\infty$. The basic analysis is done for compactly supported and rearranged data. Then we perform a series of extensions of the results to greater generality.

\medskip

\noindent $\bullet$ { \bf Compactly supported and rearranged data.}  (i) {\sl Existence and conservation of mass.} Assume to fix ideas that $u_0$ is  supported in the interval $[-R,R]$, and has total mass $M$. We may  use shifting comparison with the Barenblatt solutions with the initial masses located on either end of the support to prove bounds from above and below for the mass function of $u$, defined by $v(x,t)=\int_{-\infty} u(x,t)\,dx$, in terms of displaced versions of the mass function of the Barenblatt solution with the same mass, $V_M$ corresponding to $U_M$. We get
$$
V_M(x-R,t)\le v(x,t)\le V_M(x+R,t).
$$
Recall that $V_M(\infty,t)=M$ for every $t>0$ by mass conservation. This not only proves non-trivial existence for all times,  but also conservation of the total mass and existence of a nontrivial mass (i.\,e., $L^1$ integral) on every interval of length larger than $2R$.

We have to justify that shifting comparison applies to Barenblatt solutions, and this is done by starting the construction with compactly supported $\tilde u_{1}(x)$ and its scalings $u_{0,L}(x)=Lu_1(Lx)$.

\medskip

\noindent (ii) {\sl Positivity.} By the Aleksandrov principle, see Corollary \ref{aleks}, the function $u(x,t)$ is monotone decreasing in $x$ for $x>R$ and $t>0$. Together with the previous mass analysis this implies that $u(x,t)>0$ for $x>R$.  A similar argument happens for $x<-R$ and ensures for positivity for all $|x|>R$ and all $t>0$.

In order to establish the positivity for $|x|\le R$ at times $t>0$ we take time $t_1>0$ and move the
the origin or coordinates to a point $x_1>2R$. Setting $y=x-x_1$ we see that $u(y+x_1,t)$  is positive near the new origin and bounded below by a rearranged function $\widetilde u_0(y)$ with small support $[-R_1,R_1]$  to which the preceding result applies so that $\widetilde u(y,t)$ is positive for $|y|>R_1$. By comparison $u(x,t)$ is positive
with a uniform lower bound for $|x|\le  R$ and $t=2t_1$.

\medskip

\noindent (iii) {\sl Asymptotic behaviour as $|x|\to\infty$.}
\begin{prop} For every solution with rearranged and compactly supported   data with mass 1, we have
\begin{equation}\label{as.xbdds}
 C_1 |x|^{2s/(1+n)}\le u(x,t)\,t^{-1/(1+n)}\le C_2 |x|^{2s/(1+n)}
\end{equation}
for all $|x|\ge 2R$, and the constants do not depend on the particular $u$, nor on $t$.
\end{prop}

\noindent {\sl Proof.} The upper bound comes from the mass analysis of (i) and the monotonicity in $x$ for large $|x|$. We get
$$
x\,u(x,t)\le 2\int_{x/2}^x u(x,t)\,dx \le 2\int_{x/2}^\infty u(x,t)\,dx
$$
and by the Shifting comparison result, this mass is less than the mass of the displaced Barenblatt $U_\infty(x-R,t)$
which is proportional to
$$
(c_\infty/\gamma) \,t^{1/(1+n)}|x-R|^{-\gamma}, \quad \gamma=(2s-1-n)/(1+n).
$$
We conclude that there is a constant $c_1$ such that for all $x\ge 2R$.
$$
u(x,t)\le c_1t^{1/(1+n)}|x|^{2s/(1+n)}\,.
$$

\noindent $\bullet$ The lower bound follows from similar arguments, but now we compare with the mass of $U_\infty(x+R,t)$. We get
$$
(k-1)xu(x,t)\ge \int_{x}^{kx} u(x,t)\,dx = \int_{x}^{\infty} u(x,t)\,dx -\int_{kx}^{\infty} u(x,t)\,dx=I_1-I_2\,.
$$
The first integral $I_1$ is estimated from below  by the mass of $U_\infty(x+R,t)$ in the same interval which is accurately given by $C_2t^{1/(1+n)}|x|^{-\gamma}$, while the second is estimated from above by the mass $U_\infty(x-R,t)$ and gives
$C_1t^{1/(1+n)}|kx|^{-\gamma}$ in first approximation. Hence, for $k>C_2/C_1$ and $x$ large enough we get $u(x,t)\ge c_3t^{1/(1+n)}|x|^{2s/(1+n)}$, which implies the stated lower bound.

\medskip

\noindent {\bf Remark.} The translation of these results  for data with mass $M\ne 1$ is easy by using the mass-changing transformation $\cal T$.

\medskip

\noindent $\bullet$  {\bf Rearranged data.} If the initial data are rearranged but not compactly supported, we use approximation of the data from below with compactly supported data, so that by comparison the property of positivity follows. Conservation of mass comes from the  $L^1$ contraction property. The asymptotic lower bound in \eqref{as.xbdds} still holds, but the upper bound need not hold (it depends on the behaviour of the initial data for large $|x|$).

\medskip

\noindent $\bullet$  {\bf Continuous data.}
 In this cases there is a possible problem with the positivity of the solution, that is still not guaranteed by us. This is proved for continuous data $u_0$ by putting below some rearranged  data with compact support, $\tilde u_0$, to which the previous theory applies. Then we can apply comparison to see that $u$ has  large tails, above the minimum decay estimate of \eqref{as.xbdds}. Then conservation of mass is true for $t\ge\tau>0$.

In order to get conservation of mass since the beginning we add a small perturbation with suitable tail to $u_0$. After an easy argument this implies general conservation of mass and general minimum decay estimate.

\medskip

\noindent $\bullet$  {\bf  General data with compact support.} We attack the general case for data with compact support in the interval $I_R=[-R,R]$. We propose to use convolution with a smooth kernel to obtain a smooth approximation $u_{0\delta}$, that produces  a solution to which  the previous paragraph applies. By $ L^1 $ contraction we get conservation of mass, so that the solution must be global in time.

Due to mass conservation and the smoothing effect the mass of the solution $u$ cannot be contained in the original interval for large times. Then there is a time $T_1$  such that half the mass is outside $I_R$. By the space monotonicity away from $I_R$, we can put a displaced rearranged subsolution below, and then there are tails with at least the minimal decay rate for  $t>T_1$.  Now we use the monotonicity in time to derive the same conclusion also for all $t\ge \tau>0$.

\medskip

\noindent $\bullet$  Finally, when $u_0$ is not compactly supported, rearranged or continous, we just approximate from below with compactly supported data and pass to the monotone limit. The rest of the argument follows.

\subsection{$L^1$ continuity and initial data} We want to show that the limit solution is continuous as  an orbit $t\mapsto u(\cdot,t)$, as stated in Theorem \ref{thm.1}. At $t=0$ this means that it takes the initial data in $L^1(\re)$. We may use the fact that this holds for the approximate problems and then pass to the limit, using a Fatou argument plus conservation of mass.

The continuity at $t>0$ can be made into a stronger result, and in fact we obtain Lipschitz continuity of the orbit for all positive times. This is a consequence of  the time monotonicity and the conservation of mass. Indeed, the first implies that
for $h>0$
$$
u(x,t+h)-u(x,t)\le ((1+(h/t))^{1/(1+n)}-1) u(x,t)\le Ch\,u(x,t)\,
$$
which is uniform if $t\ge \tau>0$ and $h\le c\tau$. Since $u$ is bounded, this implies a  pointwise Lipschitz bound from above. The $L^1$ bound from below comes from conservation of mass.

\subsection{Upper limit solutions are very weak solutions}

This subsection extends to the singular case the result of papers \cite{pqrv1, pqrv2, pqrv3}. Since the approximate solutions $u_\ve$ of Problem \eqref{eq:FrUF_approx} are smooth, they satisfy the very weak formulation. The main difficulty
in passing to the limit when $\ve\to 0$ is the control of the possible growth of $\Phi_n(u)$ as $|x|\to\infty$. But this depends on having a good lower bound for $u$, and such a bound is contained in the left-hand side of \eqref{as.xbdds}. It follows that the integrals involved in the passage to the limit are uniformly absolutely integrable.

Now we can pass to the limit is the $\ve$ approximate equations written in very weak form to show that the  limit solution is indeed a very weak solution.

\subsection{Upper limit solutions form a semigroup}
We consider the maps $S_t: L^1_+(\re)\to L^1_+(\re)$ that map any initial data $u_0$ to the solution of the equation at time $t>0$, $S_t u_0= u(\cdot ,t)$. Since the approximate problems \eqref{eq:FrUF_approx} produce unique classical solutions, the semigroup property is true for them, $S^{\eps}_{t+t_1}u_0=S^{\eps}_{t}(S^{\eps}_{t_1}u_0)$. Passing to the limit we obtain the same property for our upper limit solutions arguing as follows:

Fix $u_0$ and $t_1>0$ and recall the monotone convergence
$u_\eps(x,t_1)\to u(x,t_1)$ that takes place in $L^1_{loc}(\re)$. Consider now $v_0=u(\cdot,t_1)$ as initial data for a new lap of the evolution. In order to get the upper limit solution $S_t(v_0)$ we take $\eps'>0$ and solve the approximate problem with $v_{0,\eps'}=v_0+\eps'$. If we now take $\eps'$ smaller than $\eps/2$ we get an estimate of $\|(v_{0,\eps'}- u_\eps(x,t_1))_+\|_1$ in the following way. We first note that set of points $K$ where $u(x,t_1)>\eps/2$ has measure less than $2\|u(x,t_1)\|_1/\eps$. Therefore
$$
\int_K (v_{0,\eps'}- u_\eps(x,t_1))_+\,dx =\int_K (v_{0}(x)+\eps'- u_\eps(x,t_1))_+\,dx\le |K|\eps'\,.
$$
On the other hand, since $u_\eps(x,t_1)\ge \eps$ everywhere, and $v_{0,\eps'}(x)\le \eps/2+\eps'$ on $\re\setminus K$, on that set we have $v_{0,\eps'}- u_\eps(x,t_1)\le 0$. We conclude that
$$
\|(v_{0,\eps'}- u_\eps(x,t_1))_+\|_1\le \frac{\eps'}{\eps}\|u_0\|_1\,.
$$
By the ordered contraction property, this estimate remains true during the evolution, hence
$$
\|(S^{\eps'}_t v_{0,\eps'}- u_\eps(x,t+t_1))_+\|_1\le \frac{\eps'}{\eps}\|u_0\|_1
$$
holds for all $t>0$. Let now $\eps'\to 0$ to get by the very definition of the upper limit  solution that
$$
\|(S_t v_{0}(x)- u_\eps(x,t+t_1))_+\|_1\le0.
$$
This means that $S_t v_{0}(x)\le u_\eps(x,t+t_1))$ for a.e. $ x\in \re$. In other words,
$S_t (S_{t_1} u_{0})\le S_{t+t_1}u_0$. By the conservation of mass, both functions are the same, which proves the semigroup property $S_t (S_{t_1} u_{0})= S_{t+t_1}u_0$. \qed

\subsection{Upper limit solutions are maximal solutions}

In the standard Laplacian case $s=1$ with singular diffusion, it is known that in the range of exponents $0<n<1$ where there is existence of nontrivial solutions, the limit solution is not the only possible solution defined in the whole line $x\in\re$. On the contrary, there are infinitely many other solutions with the same initial data determined by some ``flux at infinity'', as described in \cite{ERV88, RV90}. But the upper limit solution is the maximal element in that class, in fact it is the only one for which the total mass does not decrease in time.

In order to repeat the proof of maximality in a short way, we only to consider a class of solutions that admits comparison with classical supersolutions (which will be the solutions $u_\eps$ of the approximate problems). Let us call good solutions the elements of such a class.

\begin{theor}\label{th:non1ex}
 Let $u_g(x,t)$ denote a good solution of \eqref{geq.Phi} for our choice of $\Phi$ defined in an interval $0<t<T$ and having nonnegative initial data $u_0\in L^1(\re)$, and let $\overline{u}(x,t)$ the upper limit solution with same initial data. Then, $u_g\le \overline{u}$ in $\re\times (0, T)$. \end{theor}

The further exploration of the existence of such solutions falls out of the scope of this paper for reasons of space.


\section{Asymptotic behaviour of general solutions}\label{sec.asymp}

We want to prove Theorem \ref{thm.4}, which means that we want to show that the Barenblatt solutions are asymptotic attractors of the solutions with general data restricted only by the running conditions on the initial data: $u_0(x)  \ge 0$, $u_0\in L^1(\re)$. We address first the question of convergence in the $L^1$ norm, which is split into a number of cases.

\medskip

\noindent $\bullet$   We consider first  rearranged initial data. The result is just a reformulation of the proof of construction of the Barenblatt solution in Section \ref{sec.bar1d}.  The argument is well-known in the literature, see \cite{Vapme}. We just recall that the rescaled family $u_L(x,t)$ converges to $U_M$ in $L^1(\re) $ at any time $t>0$, and fix $t=1$ to get
$$
\lim_{L\to\infty}\| u_L(x,1)-U_M(x,1)\|_1\to 0
$$
We then undo the scaling, see formula \eqref{mc.scaling}, to find that
$$
\lim_{L\to\infty}\| u(x,L^{2s-1-n})-U_M(x,L^{2s-1-n})\|_1\to 0
$$
Putting $t=L^{2s-1-n}\to\infty$ we obtain the result.

\medskip

\noindent $\bullet$   The main novelty lies in establishing the result for general data that are not rearranged. We have to use the trick  introduced in paper \cite{povaz2013} with  Portilheiro for general data with compact support. We argue as follows: we  fix $t_1\gg1$, and assume the support of $u_0$ be included in  $[-R,R]$. Define
\[
    \tilde u_1(r)  := \inf_{|x|=r} u(x,t_1), \quad
    \tilde u_2(r)  := \max_{|x|=r} u(x,t_1).
\]
We easily verify that $\tilde u_1(r),\, \tilde u_2(r)$ are nonnegative and radially symmetric functions, they are nonincreasing as functions of $r$ for $r\geq R$, we have the immediate comparison
$$
\tilde u_1(r)\le u(x,t_1)\le \tilde u_2(r)\,,
$$
where $|x|=r$. We also have by the Aleksandrov reflection comparison
\[
    \tilde u_2(r)\geq \tilde u_1(r)\geq \tilde u_2(r+2R)
\]
for all $r\geq R$. It is then easy to verify that the 1-$d$ mass of $\tilde u_2(r) - \tilde u_1(r)$ is less than
$CRt_1^{-1/(n+1)}$, which can  be made very small.

We restart the evolution at time $t_1$ and get radially symmetric solutions $  \tilde u_1(x,t),   \tilde u_2(x,t)$ with initial data $  \tilde u_1(r),   \tilde u_2(r)$ resp., and we also have the original $u(x,t+t_1)$ (now displaced in time) that stays between them.
The asymptotic behaviour says that $\tilde u_1(x,t)$ converges to the Barenblatt $U_{M_1}$, and $\tilde u_1(x,t)$ converges to the Barenblatt $U_{M_2}$. Moreover, the masses satisfy $M_1\le M\le M_2$ and $M_2-M_1\le \eps$. The asymptotic formula for convergence of $u(\cdot,t)$ to $U_M(\cdot,t)$ follows easily.

\medskip

\noindent $\bullet$   If $u_0$ does not have compact support we use approximation and $L^1$ contraction.

\medskip

\noindent $\bullet$   The estimate in the $L^p$ norms is just an interpolation between the convergence result for $p=1$ just proved, and the  $L^\infty$ bound of the form $u(x,t)\le C\,t^{-\alpha}$ with $\alpha=1/(2s-1-n)$.\qed

\medskip

This result allows to extend the uniqueness of the fundamental solutions as follows

\begin{theor} Every upper limit solution for positive times that that is self-similar and integrable is a Barenblatt solution $U_M(x,t)$ for some $M>0$. It can be obtained by rescaling from any integral and nonnegative initial data with mass $M$.
\end{theor}

\noindent {\sl Proof.} Let $u_1$ be that solution and $M$ its mass. Since it is self-similar  we have
$$
\|u_1(\cdot,t_1)-U_M((\cdot,t_1)\|_1=\|u_1(\cdot,t_2)-U_M((\cdot,t_2)\|_1
$$
Fix now $t_1>0$ and let $t_2\to \infty$. Applying Theorem \ref{thm.4} the right-hand side goes to zero. Hence, $u_1 \equiv U_M$. The second assertion is easier. \qed

\section{Proof of the Shifting comparison lemma}\label{sec.shifting}

In this section we will give a complete proof of Theorem \ref{thm.shift}.
 Note that this result applies to the approximate problems that have nonsingular functions $\Phi$. The passage to the limit allows to apply it to our upper limit solutions.

\subsection{Elliptic problem. Extended problem}

The implicit time discretization scheme \cite{BeTh72, CL71} directly connects the analysis of the parabolic equation \eqref{geq.Phi} to solving a sequence of elliptic equations of the form
\begin{equation} \label{eq.ell1}
\left(  -\Delta\right)^{\sigma/2}v+  B(v)=f\left(  x\right)   \qquad x\in \re,%
\end{equation}
where  $\sigma\in(0,2)$ and $f$ is an integrable function defined in $\re$. We assume that the nonlinearity is given by a function $B:\R_{+}\rightarrow\R_{+}$ which is  smooth and monotone increasing with $B(0)=0$ and $B'(v) >0$. It is not essential to consider negative values for our main results, but the general theory can be done in that greater generality. We are using here nonnegative data and solutions. In the parabolic application $B=\Phi^{-1}$, see \cite{Vapme, VazVol1}. Then we need to prove the following result

\begin{theor}
Let us consider two functions $f_1, f_2\in L^1(\re)$ with the following properties:

(i) They are  nonnegative,  and rearranged  around their points of maximum $x_1$ and $x_2$ resp., with $x_1<x_2$. BY this mean that $f_i(x-x_i)$ is rearranged.

(ii) We assume moreover that $\int_{-\infty}^x f_2\,dx\le \int_{-\infty}^x f_1\,dx$

(iii)  the mass is the same, $\int f_1(x)\,dx=\int f_2(x)\,dx=M>0$.

\noindent Besides, we assume that $B$ is convex and defined on $\re_+$ with $B(0)=0$ and $B'(0)>0$. Then,  the following comparison inequalities hold
\begin{equation}
\int_{-\infty}^x v_2(x)\,dx\le \int_{-\infty}^x v_1(x)\,dx, \qquad \int_{-\infty}^x B(v_2(x))\,dx\le \int_{-\infty}^x\, B(v_1(x))\,dx
\end{equation}
for every $x\in\re.$
\end{theor}

\noindent {\sl Proof.} (a) It will be convenient to formulate the elliptic problem by using a proper extension problem, which  is defined as the trace of a properly defined Dirichlet-Neumann problem as follows. If $w$ is a weak solution to the local problem
\begin{equation}\label{eq.3}%
\left\{
\begin{array}
[c]{lll}%
\operatorname{div}_{x,y}\left(  y^{1-\sigma}\nabla w\right)  =0 &  & \text{in} \ Q_+,\\[10pt]
\displaystyle{-\frac{1}{\kappa_{\sigma}}\lim_{y\rightarrow0^{+}}y^{1-\sigma}\,\dfrac{\partial w}{\partial y}(x,y)}+\,B(w(x,0))=f(x)   & & \text{for } \ x\in \re
\end{array}
\right.
\end{equation}
where $Q_+:=\re\times\left(  0,+\infty\right)  $ is the upper half-plane and $\kappa_{\sigma}$ is the constant  that is not important in what follows. See \cite{CS2007}which is the main reference in the issue of this extension. We can define again a suitable meaning of weak solution in terms of this extended problem. The functional setting is perfectly explained in \cite{VazVol1}, Section 3.1. There is a solution to problem (\ref{eq.3}). Then the trace of $w$ over $\re\times \{0\}$, $\text{Tr}_{\re}(w)=w(\cdot,0)=:v$ is said a solution to problem (\ref{eq.ell1}).
Using the change of variables $z=Cy^{\sigma}$ for a convenient constant $c>0$, the problem  can also be written as
\begin{equation}
\left\{
\begin{array}
[c]{lll}%
c\,z^{\nu}\dfrac{\partial^{2}w}{\partial z^{2}}+\dfrac{\partial^{2}w}{\partial x^{2}}=0 &  & \text{ in
} \ Q_+\\[10pt]
-\dfrac{\partial w}{\partial z}\left(  x,0\right)  =%
f(x)-B(w(x,0)) &  & \text{for } \ x\in \re\,.
\end{array}
\right.  \label{eq.23}%
\end{equation}
 We use the extension formulation and write $w_i(x,0)=v_i(x)$ and $u_i=B(v_i)$, $i=1,2$.

(b) Note that the properties (i) of the $f_i$'s hold also  for the solutions $v_1$ and $v_2$. Conservation of mass applies so that $B(v_1)$ and $B(v_2)$ have the same mass $M$. Due to the properties of $B$ the functions $v_1$ and $v_2$ are also integrable, though their masses are not controlled. By the properties of the extension to the upper half-plane, the functions $w_2(z,z)$ and $w_1(x,z)$ are also integrable in $x$ for every fixed $z>0$.

(c) Using a similar strategy to the symmetrization proof in paper \cite{VazVol1}, we introduce the function
\begin{equation}
{Z}(s,z)=\int_{-\infty}^{s}(w_2(\tau,z)-w_1(\tau,z))d\tau\,.
\end{equation}
Then, it is clear that
\begin{equation}
c\,z^{\nu}{Z}_{zz}+ {Z}_{xx}=0 \label{symineq}
\end{equation}
and
\begin{equation}
Z(-\infty,z)=0, \qquad {Z}(\infty,z)=0\,.\label{boundcond}
\end{equation}
A crucial point in our arguments below is played by the derivative of $Z$ with respect to $z$. Due to the boundary conditions contained in \eqref{eq.23}, we have
\begin{equation}\label{Z_yboundary.formula}
{Z}_{z}(x,0)\geq \int_{-\infty}^{x} (B(w_2(\tau,0))-B(w_1(\tau,0))\, d\tau
\end{equation}
Recall that $w_i(\tau,0)=v_i(\tau)$. Observe also that the function
$$
Y(x,0)=\int_{-\infty}^{x}B(w_2(\tau,0))-B(w_1\tau,0))\, d\tau
$$
has the same points of maximum or minimum and the same regions of monotonicity than ${Z}(x,0)$.

\noindent (d) Then we argue as follows. Due to the maximum principle and the boundary conditions \eqref{boundcond}, a positive maximum of $Z$ can be achieved only  on the line $\left\{z=0\right\}$. On the other hand, in the interval $I= \{x: x_1< x< x_2\}$ we know that $\partial_x v_1\le 0$ and $\partial_x v_2\ge 0$ hence in this interval $I$
$$
Y_{xx}=\partial_x B(v_2)-\partial_x B(v_1)\ge 0
$$
and the maximum of $Y$ must lie outside of $I$, hence the same happens for $Z$. Suppose the maximum  of $Z$ happens at $(x_{0},0)$ with $x_0\ge x_2$. We must also have $Z_z(x_0,0)<0$ by Hopf's maximum principle, and by \eqref{Z_yboundary.formula}, this
leads to $Y(x_0,0)<0$. But for $x>x_0$
$$
\begin{array}{l}
\displaystyle Y(x,0)-Y(x_0,0)=\int_{x_0}^x [B(v_2(\tau))- B(v_1(\tau))]\, d\tau\\
[6pt]
\displaystyle \le \int_{s_0}^s
B'(v_2(\tau,0))(v_2(\tau)-v_1(\tau))\, d\tau.
\end{array}
$$
Here, we have used the convexity of $B$ so that $B'$ is an increasing real function and
$$
B(v_2(\tau))-B(v_1(\tau))\le B'(v_2(\tau))\,(v_2(\tau))-v_1(\tau))\,.
$$
After integration by parts in the expression for the increment of $Y$, we get
$$
\begin{array}{l}
\displaystyle Y(x,0)-Y(x_0,0) \le \left[B'(v_2(\tau))(Z(\tau,0)-Z(x_0,0))\right]_{x_0}^{x} -\\
[6pt]
\displaystyle \int_{x_0}^x B''(v_2(\tau))v_{2,x}(\tau)(Z(\tau,0)-Z(x_0,0))d\tau.
\end{array}
$$
Since $Z$ has a maximum at $x_0$ and $B'$ is positive, the first term in the RHS is non-positive. As for the second, we have: $B''>0$, $v_{2,x}<0$, and  $Z(x,0)-Z(x_0,0)\le 0$,  hence the last term is also nonpositive. We conclude that $Y(x,0)\le Y(x_0,0)<0$ for all $x>x_0$. This is a contradiction, because by the
conservation of mass property at plus infinity we have $Y(\infty,0)=0$. Therefore, there is no positive maximum for $Z$ on this side.

(ii) Similar argument on  the other side, $x<x_1$ reversing the roles of $v_1$ and $v_2$ and the direction of integration, that starts now at $+\infty$. We conclude that  $Z(x,0)\le 0$ everywhere.

(iii) Once we have $Z(x,0)\le 0$ we also want to prove that  $Y(x,0)\le 0$.  We may use Lemma \ref{lemma1} below (a well-known result), taking advantage of the convexity of $B$  and choosing any convex, increasing
function $\Phi:[0,\infty)\rightarrow[0,\infty)$. This ends the proof of the  comparison theorem in this case. \qed

\begin{lemma}\label{lemma1}
Let $f,g\in L^{1}(\Omega)$ be two rearranged functions on a ball $\Omega=B_{R}(0)$. Then $f\prec g$ if and only if for every convex
nondecreasing
function
$\Phi:[0,\infty)\rightarrow [0,\infty)$ with $\Phi(0)=0$ we have
\begin{equation}
\int_{\Omega}\Phi(f(x))\,dx\leq \int_{\Omega}\Phi(g(x))\,dx.
\end{equation}
This result still holds if $R=\infty$ and $f,g\in L^{1}_{loc}(\R^{N})$ with $g\rightarrow0$ as $|x|\rightarrow\infty$.
\end{lemma}

\medskip

\noindent {\bf Remark.}  We have imposed severe conditions on the shape of $u_{01}, u_{02}$, and we have required $\Phi$ to be concave. Are these conditions necessary? They are not for the same result with standard Laplacian instead of fractional Laplacian.


\section{A preview of logarithmic diffusion with $s=1/2$}\label{sec.logdiff}

This is a kind of exceptional case in the parameter diagram. The study  is a bit different from the previous analysis, and much of our intuition comes from a similar special case that happens  for standard diffusion $s=1$, in dimension $N=2$ with logarithmic diffusion $n=0$. An explicit solution exists then and it was used in \cite{Vaz92}. The formula is
\begin{equation}
U(x,t)=\frac{8(T-t)}{(1+|x|^2)^2}, \quad x\in \re^2, \ 0<t<T.
\end{equation}
Note that $U$ is positive only in the time interval $0<t<T$. As a very weak solution, it can be defined for all $t>0$ and  it vanishes identically for $t\ge T$. Another important property is  the total mass decay with a constant rate
\begin{equation}
\frac{d}{dt}\int_{\re^2} U(x,t)\,dx= -8\pi.
\end{equation}
But, as shown by a number of studies, this rate does not correspond the limit solutions, which are characterized by the rule $dM(t)/dt=-4\pi$. A  gap of non-uniqueness opens up and this is carefully described in \cite{JLVSmoothing}, where related literature can be found

\subsection{Existence of a positive solution}

 Luckily, in our one-dimensional case  there also exists an explicit solution  of the evolution equation
\begin{equation}
\label{eq:logonehalf} \partial_t u + \Doh(\log u) = 0\,\,\,\,\hbox{ in } \R\times (0,T)
\end{equation}
with a smooth initial condition $u(x,0)\in L^1(\R)$. It is given by the formula
\begin{equation}\label{explicit.12}
U(x,t)=\frac{2\lambda\,(T-t)}{\lambda^2+|x|^2} \,\,\,\hbox{ in } \R\times (0,+\infty),
\end{equation}
with any $\lambda>0$, so it is indeed a whole family of solutions related by scaling.
We see that $U$ has the separate-variable type as in the previous example, and $U(\cdot, t)$ is in $L^1(\R)$ for any $0\le t< T$. It is very peculiar that the solution becomes identically zero in finite time. This is the so-called {\sl finite-time extinction} phenomenon which is typical of some ranges of fast diffusion, see \cite{JLVSmoothing} for standard diffusion and \cite{pqrv2, KL} for fractional diffusion.

In order to prove that  \eqref{explicit.12} is a classical solution of the equation we write it in the form $U(x,t)=(T-t)F(x)$ and we need to find an integrable  profile $F>0$ such that
\begin{equation}\label{ell.log}
\Doh F(x)=F(x).
\end{equation}
It only remains to check that $F(x)=2/(1+|x|^2) $ is a solution of this nonlinear elliptic equation. This is only a technical calculus result.

Observe that the initial mass is in all cases \eqref{explicit.12}
$$
M_0=\|U_0\|_1=2T\int \frac{dx}{1+x^2}=2\pi T,
$$
so that in the extinction time is given by $T=\|U_0\|_1/{2\pi}.$ Accordingly, the mass decay rule is $M'(t)=-2\pi$, which is related to what happens for $N=2$ and standard diffusion (see above).

In the same way as in the previous analysis of the case $s>1/2$, we can use this example and standard comparison to prove that for initial data  $u_0(x)\ge c/(1+|x|^2)$ with $c>0$, the limit solution is nontrivial, more precisely positive for $0<t<T$ with $T=c/2$. We do not expect upper limit solutions to conserve mass since this does happen for the special case of standard diffusion $s=1$, $N=2$, as described in \cite{rev97, VER96, JLVSmoothing}.

\section{Comments, extensions, and open problems}\label{sec.ceop}

\ $\bullet$ We have given preference in the paper to the treatment of the more singular case of exponents $n>0$. However, the logarithmic equation is also covered and the main results are proved to be true, as special case $n=0$. This is shown not only at the formal level, but also technical details are given as needed.  There is another variant of the limit case $n=0$, the sign diffusion $u_t=\Delta \mbox{\rm sign\,}(u)$, that is related to the total variation flow and was treated in \cite{BoFi2012}. We are not covering the fractional version of this variant.

\medskip

$\bullet$ The actual behaviour of the solutions of logarithmic diffusion on the line with exponent $s=1/2$ is a question to be investigated, and this will done in a separate work. For the non-singular equation $\partial_t u + (-\Delta)^s(1+\log u) = 0$ the study was done in \cite{pqrv3}, see also \cite{VPQR13}.

\medskip

$\bullet$  The existence of non-maximal solutions, that are not upper limit solutions and do not conserve mass, is an interesting open problem. For the standard Laplacian, it was solved in \cite{ERV88, RV90} where it is proved that there are infinitely many solutions for every integrable initial data and they are determined by some flux conditions at infinity. A large related literature has developed, see e.\,g. \cite{Hui2004, Hsu2002, RV90, RV92}.

\medskip

$\bullet$ Equation \eqref{eq.mass} is a kind of fractional $p$-Laplacian equation for the mass function.
It would be interesting to perform a study of its properties and applicability.

\medskip

$\bullet$ Elliptic problems and the Crandall-Liggett approach. A natural way to solve the evolution equation is by implicit discretization in time, a method that became basic in the early studies of the Porous Medium Equation, see the original paper {CL71} or \cite[Chapter 10]{Vapme}. In the present situation, it means that we have to solve elliptic problems of the form $(-\Delta)^s \Phi(u)+ u=f$. Note that when $\Phi(u)=\log(u)$, this equation takes the suggestive form $(-\Delta)^s v+ e^v=f$. Non-existence results are carefully described in \cite{BSV15} for the range of parameters $2s<n+1$. The existence theory in our parameter range is not difficult. We refrain from further details for reasons of space, but see \cite[Section 9]{BSV15}.

\medskip

$\bullet$ Dirichlet problem in bounded domains. In the case of the standard Laplacian, $s=1$, no nontrivial solutions exist for the Dirichlet problem with zero boundary data. The question is open for our equations.

\medskip

$\bullet$  Problems with more general nonlinearities $\Phi$. Symmetrization can be applied  in some cases using the results of \cite{VazVol2} to compare a general $\Phi$ with the power cases we deal with. It can give the starting results on nontrivial existence. Non-existence results in that direction are explained in \cite{BSV15}. Of course, in the simple case where $\Phi(u)=-c\,u^{-n}$ (or $\Phi(u)=c\,\log(u)$) with some $c>0$, the constant may be absorbed into the time variable and we need not make any changes to the theory.

\medskip

$\bullet$  Problems with other classes of initial data are worth studying. A simple example are the constant solutions do not belong to our class of integrable upper limit solutions. It is not difficult to construct solutions for initial data in $L^p(\re)$ with $p>1$ using the smoothing effect and approximation by integrable data. We can  consider in this way initial data that are not integrable with different decay rates at infinity. It is even possible to consider data that grow at infinity.

\medskip

$\bullet$  A different direction is establishing existence of upper limit solutions with nonnegative Radon measures as initial data. An example for that extension are the Barenblatt solutions that we have just constructed. Since the basic process for general measures is easy following the indications of \cite[Section 4]{VazBar}, we leave it to the interested reader.

\vskip .8cm

\noindent{\bf Acknowledgments. } Work partially supported by Spanish Project MTM2011-24696. Part of this work was done  at the Isaac Newton Institute, Cambridge, UK,  during the program Free Boundary Problems and Related Topics held in the spring of 2014. The author thanks A. de Pablo for information about the explicit solution \eqref{ell.log}.

\vskip 1cm



\

\

\noindent {\sc Address:}

\medskip

\noindent {\sc Juan Luis V{\'a}zquez}\newline
Departamento de Matem\'{a}ticas, Universidad Aut\'{o}noma de Madrid, \\
28049 Madrid, Spain. \ e-mail: {\tt juanluis.vazquez@uam.es}

\vskip 1cm


{\small
\addcontentsline{toc}{section}{~~~References}
\begin{thebibliography}{10}

\bibitem{Alek60} {\sc A.~D. Aleksandrov.} {\sl Certain estimates for the Dirichlet
problem}, {Soviet Math. Dokl.} {\bf 1} (1960), 1151--1154.

\bibitem{ArBk86} {\sc D.~G. Aronson.} {\sl The porous medium equation}, CIME Lectures, In {\sl ``Some problems in nonlinear diffusion''} (K. Kirchg\"assner H. Amann, N. Bazely, editors), Lecture Notes in Mathematics {\bf 1224}, Springer-Verlag, New York, 1986.

\bibitem{Bandle} {\sc C. Bandle}.
{\sl Isoperimetric inequalities and applications.}
Monographs and Studies in Mathematics, 7. Pitman
(Advanced Publishing Program), Boston, Mass.-London, 1980.

\bibitem{BeTh72} {\sc P. B\'{e}nilan}.
{\it Equations d'\'evolution dans un espace de Banach quelconque
et applications}, Ph. D. Thesis, Univ. Orsay, 1972 (in French).

\bibitem{BenCr} {\sc P. B\'enilan, M. G. Crandall.}
{\sl Regularizing effects of homogeneous evolution equations.}
Contributions to Analysis and Geometry, (suppl.  to Amer. Jour.
Math.) Johns Hopkins Univ. Press, Baltimore, Md., 1981. Pp.
23-39.

\bibitem{BoFi2012} {\sc M. Bonforte, A. Figalli.} {\sl Total variation flow and sign fast diffusion in one dimension}, J. Differ. Eqns. {\bf 252} (2012), no. 8, 4455--4480.

\bibitem{BSV15} {\sc M. Bonforte, A. Segatti, J.~L. V{\'a}zquez}. {\sl Non-existence and instantaneous extinction of solutions for singular nonlinear fractional diffusion equations}, Preprint.

\bibitem{BV2012} {\sc M. Bonforte,  J. L. V\'azquez.}
{\sl Quantitative Local and Global  A Priori Estimates for Fractional Nonlinear Diffusion Equations,} Advances in Math. {\bf 250} (2014) 242--284.

\bibitem{BV2015} {\sc M. Bonforte,  J. L. V\'azquez.} {\sl A Priori Estimates  for Fractional Nonlinear  Degenerate Diffusion Equations  on bounded domains}. To appear in Arch. Rat. Mech. Anal. (2015). \texttt{http://arxiv.org/abs/1311.6997}

\bibitem{CS2007} {\sc L.A. Caffarelli, L. Silvestre. } {\sl An extension problem related to the fractional Laplacian,} {Comm. Partial Differential Equations,} \textbf{32} (2007), 1245--1260.

\bibitem{CVW87} {\sc  L. A. Caffarelli, J. L. V\'azquez, N. I. Wolanski.}
{\sl Lipschitz continuity of solutions and interfaces of the $N$-dimensional porous medium equation},
{Indiana Univ.\ Math.\ J.} {\bf 36} (1987), no. 2, 373--401.

\bibitem{ChaVaz2002} {\sc  E. Chasseigne,  J.~L. V\'azquez.} {\sl Theory of
Extended Solutions for Fast Diffusion Equations in Optimal Classes
of Data. Radiation from Singularities},
 Archive Rat. Mech. Anal. {\bf 164} (2002), 133--187.

\bibitem{CL71} {\sc M.~G. Crandall, T.M. Liggett.}
{\sl Generation of semi-groups of nonlinear transformations on general
Banach spaces},  { Amer. J. Math.} {\bf 93} (1971) 265--298.

\bibitem{DDP1994} {\sc P. Daskalopoulos,  M. A. del Pino}.
{\sl On fast diffusion nonlinear heat equations and a related singular elliptic problem}, Indiana Univ. Math. J. 43 (1994), no. 2, 703--728.

\bibitem{DDP1995} {\sc P. Daskalopoulos, M. A. del Pino}.
 {\sl  On a singular diffusion equation,} Comm. Anal. Geom. 3 (1995), no. 3-4, 523–542.

\bibitem{DDP1997} {\sc P. Daskalopoulos,  M. A. del Pino}.
 {\sl On nonlinear parabolic equations of very fast diffusion},
 Arch. Rational Mech. Anal. 137 (1997), no. 4, 363--380.

\bibitem{DaskKenig} {\sc P. Daskalopoulos, C. E. Kenig}. {\sl ``Degenerate diffusions.  Initial value problems and local regularity theory''}, EMS Tracts in Mathematics, 1. European Mathematical Society (EMS), Z\"urich, 2007.

 \bibitem{pqrv1} {\sc A. De Pablo, F. Quir\'os, A. Rodr\'{\i}guez, J. L. V\'azquez.}
{\sl A fractional porous medium equation}, Advances in Mathematics
226 (2011), no. 2, 1378--1409.

\bibitem{pqrv2} {\sc A. De Pablo, F. Quir\'os, A. Rodr\'{\i}guez, J. L. V\'azquez.}
{\sl A general fractional porous medium equation}.
Comm. Pure Appl. Math. {\bf 65} (2012), no.~9, 1242--1284.

\bibitem{pqrv3} {\sc  A. de Pablo, F. Quir\'os, A. Rodr\'{\i}guez, and J.~L. V\'azquez}.
{\sl Classical solutions for a logarithmic fractional diffusion equation}.
Journal de Math. Pures Appliqu\'ees, to appear.

\bibitem{ERV88} {\sc J. R. Esteban, A. Rodr\'{\i}guez, J. L. V\'azquez.}
{\sl A nonlinear heat equation with singular diffusivity,}
 {Comm. Partial Diff. Eqs.} {\bf 13} (1988), 985--1039.

\bibitem{Herrero89} {\sc  M. A. Herrero.}
{\sl A limit case in nonlinear diffusion}
Nonlinear Anal. {\bf 13} (1989), no. 6, 611--628.

\bibitem{Hui2004} {\sc K. M. Hui}.
{\sl Existence of solutions of the very fast diffusion equation},
Nonlinear Anal. {\bf 58 } (2004), no. 1-2, 75--101.

\bibitem{Hsu2002} {\sc S.-Y. Hsu}.
{\sl Dynamics near extinction time of a singular diffusion equation},
Math. Ann. {\sl 323} (2002), no. 2, 281--318.

\bibitem{Igb09} {\sc N. Igbida}.
{\sl From fast to very fast diffusion in the nonlinear heat equation},
 Trans. Amer. Math. Soc. {\bf 361} (2009), no. 10, 5089--5109.

\bibitem{KL}{\sc S.~Kim, K.~A. Lee}.
{\sl H\"older estimates for singular nonlocal parabolic equations}.
Journal of Functional Analysis, {\bf 261}, (2011) 3482--3518.


\bibitem{povaz2013} {\sc   M. Portilheiro, J. L. V\'azquez.}
{\sl   Degenerate homogeneous parabolic equations  associated with the
infinity-Laplacian}. Calc. Var. PDE. {\bf 46} (2013), no. 3-4, 705--724.

\bibitem{RV90} {\sc A. Rodr\'{\i}guez, J. L. V\'azquez}.
{\sl A well-posed problem in singular Fickian diffusion}, {Archive
Rat. Mech. Anal.} {\bf  110}, 2 (1990), 141--163.

\bibitem{RV92} {\sc A. Rodr\'{\i}guez, J. L. V\'azquez}.
{\sl Maximal solutions of singular diffusion equations with general initial data},
Nonlinear Diffusion Equations and their Equilibrium States, 3, 7 Birkhäuser Verlag, Boston (1992) p. 471–484.

\bibitem{rev97} {\sc  A. Rodr\'{\i}guez, J.~L. V\'azquez, J.~R. Esteban}.
{\sl The maximal solution of the logarithmic fast diffusion equation in two space dimensions}, Adv. Differential Equations {\bf 2 } (1997), no. 6, 867--894.

\bibitem{Rosen95} {\sc  Ph. Rosenau.}
{\sl Fast and superfast diffusion processes}, Physical Rev. Let. {\bf 74}, 7 (1995), 1056--1059.


\bibitem {Ser71} {\sc J. Serrin}. {\sl A symmetry problem in potential theory}, {Arch.
Rat. Mech. Anal.} {\bf 43} (1971), 304-318.

\bibitem{Vsym82} {\sc J.~L. V\'azquez}.
{\sl Sym\'etrisation pour $u_t=\Delta\varphi(u)$ et applications,}
C. R. Acad. Sc. Paris {\bf  295} (1982), pp. 71--74.

\bibitem{VTr1} {\sc J.~L. V\'azquez,}
{\sl  Asymptotic behaviour and propagation properties of the
one-dimensional flow of gas in a porous medium},
{Trans. Amer. Math. Soc.} {\bf 277} (1983), pp. 507--527.


\bibitem{Vaz92}
{\sc J.~L.~V\'azquez}.
{\sl Nonexistence of solutions for nonlinear heat equations of fast-diffusion type},
J. Math. Pures. Appl. {\bf 71} (1992), 503--526.


\bibitem{VazGaeta} { \sc  J. L. V{\'a}zquez.} {\sl The Porous Medium Equation. New
contractivity results}. In    Progress in Nonlinear Differential Equations and Their
Applications, {\bf 63} (205)  (Volume in honor of H. Brezis), pp. 433--451.

\bibitem{VANS05} {\sc J.~L. V{\'a}zquez}.
{\sl Symmetrization and Mass Comparison for
Degenerate Nonlinear Parabolic and  related Elliptic Equations},
Advances in Nonlinear Studies, {\bf 5} (2005), 87--131.

\bibitem{Vapme} {\sc J.~L. V\'azquez}.
 {``The Porous Medium Equation. Mathematical Theory''},
 Oxford Mathematical Monographs. The Clarendon Press, Oxford University Press, Oxford  (2007).


\bibitem{JLVSmoothing}{\sc J.~L. V\'azquez}.
{ \lq\lq Smoothing And Decay Estimates
For Nonlinear Diffusion Equations. Equations Of Porous Medium
Type''}, Oxford Lecture Series in Mathematics and its Applications,
33. Oxford University Press, Oxford, 2006.

\bibitem{VazBar}{\sc J.~L. V\'azquez}. {\sl Barenblatt  solutions and asymptotic behaviour for a  nonlinear fractional heat equation of porous medium type. }   J. Eur. Math. Soc. {\bf 16} (2014), 769--803. MR3191976. In  arXiv:1205.6332.

\bibitem{VER96} {\sc J.~L. Vazquez, J. R. Esteban, A. Rodr\'{\i}guez}.
{\sl The fast diffusion equation with logarithmic nonlinearity
and the evolution of conformal metrics in the plane},
{\em  Advances  Diff. Eqns} {\bf 1}, 1 (1996), 21--50.

\bibitem{VPQR13}
{\sc  J.~L.~V\'azquez,  A. de Pablo,  F. Quir\'os, A. Rodr\'iguez.}
{\sl Classical solutions and higher regularity for nonlinear fractional diffusion
equations},
preprint ArXiv 1311.7427 (2013).


\bibitem{VazVol1} {\sc J.~L. V\'{a}zquez, B.~Volzone}.
{\sl Symmetrization for Linear and Nonlinear Fractional  Parabolic Equations of Porous Medium Type,} J. Math. Pures Appl. (9) {\bf 101} (2014),  no. 5, 553--582. 

\bibitem{VazVol2} {\sc J.~L. V\'{a}zquez, B.~Volzone}.
{\sl   Optimal estimates for Fractional  Fast diffusion equations}.
J. Math. Pures Appl. {\bf 103} (2015), 535--556.

\bibitem{Villani2003} {\sc  C. Villani.} {``Topics in Optimal Transportation''},  Graduate Studies in Mathematics 58, American Mathematical Society, Providence (2003).

\end{thebibliography}
}
\end{document}